\documentclass[11pt,a4paper]{amsart}
\usepackage{a4wide}
\usepackage{amsthm, amsmath, amsfonts, amssymb,graphicx}
\usepackage{mathrsfs,url}

\date{\today}
\subjclass[2010]{primary 20B27, 05C55, 05C05, 08A35, 03C40; secondary 08A70}
\renewcommand{\perp}{{\lVert}}
\renewcommand{\iff}{if and only if }
\newcommand{\dff}{if }

\theoremstyle{plain}
\newtheorem{thm}{Theorem}[section]
\newtheorem{theorem}[thm]{Theorem}
\newtheorem{lem}[thm]{Lemma}
\newtheorem{lemma}[thm]{Lemma}
\newtheorem{prop}[thm]{Proposition}
\newtheorem{proposition}[thm]{Proposition}
\newtheorem{cor}[thm]{Corollary}
\newtheorem{corollary}[thm]{Corollary}
\theoremstyle{definition}
\newtheorem{definition}[thm]{Definition}
\newtheorem{defn}[thm]{Definition}

\author[M. Bodirsky]
{Manuel Bodirsky}
    \address{Institut f\"{u}r Algebra\\TU Dresden
    \\01062 Dresden\\Germany}
    \email{Manuel.Bodirsky@tu-dresden.de}
   \urladdr{http://www.math.tu-dresden.de/~bodirsky/}
    \thanks{The first and fourth author have received funding from the European Research Council under the European Community's Seventh Framework Programme (FP7/2007-2013 Grant Agreement no. 257039). 
    The first author also received funding from the German Science Foundation (DFG, project number 622397). 
    The second author was supported by the European Community's Marie Curie Initial Training Network in Mathematical Logic - MALOA - From MAthematical LOgic to Applications, PITN GA-2009-238381 during a short term early stage research fellowship held at the \'{E}quipe de Logique Math\'{e}matique, Universit\'{e} Diderot - Paris 7, and by the German Research Foundation grant SFB 878. The third author has
been funded through projects I836-N23 and P27600 of the  Austrian
Science Fund (FWF). The fourth author was supported by the Hungarian Scientific Research Fund (OTKA) grant no.~K109185.}

\author[D. Bradley-Williams]
{David Bradley-Williams}
    \address{
    Mathematisches Institut der Heinrich-Heine-UniversitŠt,
Universit\"atsstr.\ 1,
40225 D\"usseldorf, Germany.}
\email{bwilliams.david@googlemail.com}
\urladdr{http://www1.maths.leeds.ac.uk/pure/logic/homepages/bwilliams.html}

\author[M. Pinsker]
{Michael Pinsker}
	\address{Department of Algebra, MFF UK, Sokolovska 83, 186 00 Praha 8, Czech Republic}        
	\email{marula@gmx.at}
    \urladdr{http://dmg.tuwien.ac.at/pinsker/}

\author[A.~Pongr\'{a}cz]
{Andr\'{a}s Pongr\'{a}cz}
    \address{Department of Algebra and Number Theory\\
    University of Debrecen\\
    4032 Debrecen, Egyetem square 1\\
    Hungary}
\email{pongracz.andras@science.unideb.hu}

\title{The Universal Homogeneous Binary Tree}

\newcommand{\mS}{\mathbb S}
\newcommand{\mL}{\mathbb L}
\newcommand{\mQ}{\mathbb Q}
\DeclareMathOperator{\Betw}{Betw}
\DeclareMathOperator{\Cyc}{Cyc}
\DeclareMathOperator{\Sep}{Sep}
\DeclareMathOperator{\Csp}{CSP}

\newcommand{\To}{\rightarrow}

\DeclareMathOperator{\Aut}{Aut}

\newcommand{\rest}{{\upharpoonright}}
\newcommand{\mix}[2]{\binom{#1}{#2}}
\newcommand{\st}{\mathbb S_2}
\newcommand{\ignore}[1]{}

\DeclareMathOperator{\Sym}{Sym}
\DeclareMathOperator{\Emb}{Emb}
\DeclareMathOperator{\End}{End}

\begin{document}

\begin{abstract}
A partial order is called \emph{semilinear} \dff
the upper bounds of each element are linearly ordered and any two elements have a common upper bound. 
There exists, up to isomorphism, a unique 
countable existentially closed semilinear order, which we denote by 
$(\mS_2;\leq)$. 
We study the \emph{reducts} of $(\mS_2;\leq)$, that is, the relational structures with domain $\mS_2$, all of whose relations are first-order definable in $(\mS_2;\leq)$.
Our main result is a classification of the model-complete cores of the reducts
of $\mS_2$. From this, we also obtain
a classification of reducts up to first-order interdefinability, which is equivalent to a classification of all 
subgroups of the full symmetric group
on $\mS_2$ that contain the automorphism group of $(\mS_2;\leq)$
and are closed 
with respect to the pointwise convergence topology. 

\end{abstract}

\maketitle
\section{Introduction}
A partial order $(P;\leq)$ is called 
\emph{semilinear} \dff 
for all $a,b \in P$ there exists $c \in P$ such that
$a \leq c$ and $b \leq c$, and 
for every $a \in P$ the set $\{ b \in P : a \leq b\}$ is linearly ordered, that is, contains no incomparable pair of elements. 
Finite semilinear orders are closely related to rooted 
trees: the transitive closure of a rooted tree (viewed as a directed graph with the edges
oriented towards the root) is a semilinear order,
and the transitive reduction of any finite semilinear order is a rooted tree.  

It follows from basic facts in model theory (e.g.~Theorem 8.2.3.~in~\cite{Hodges})
that there exists a countable semilinear order 
$(\mS_2;\leq)$ which is \emph{existentially closed} in the class of all countable semilinear orders,
that is, for every embedding $e$ 
of $(\mS_2;\leq)$ into a countable
semilinear order $(P;\leq)$, 
every existential formula $\phi(x_1,\dots,x_n)$,
and all $p_1,\dots,p_n \in \mS_2$ such that $\phi(e(p_1),\dots,e(p_n))$ holds in $(P;\leq)$ we have
that $\phi(p_1,\dots,p_n)$ holds in $(\mS_2;\leq)$.  
We write $x < y$ for $(x \leq y \wedge x \neq y)$ and $x \perp y$ for $\neg (x \leq y) \wedge \neg (y \leq x)$, that is, for incomparability with respect to $\leq$. Clearly, $(\mS_2;\leq)$ is 
\begin{itemize}
\item \emph{dense}: for all $x,y \in \mS_2$
such that $x < y$ there exists $z \in \mS_2$ such that $x < z < y$; 
\item \emph{unbounded}: for every
$x \in \mS_2$ there are $y,z \in \mS_2$ such that $y < x < z$;
\item \emph{binary branching}: 
(a) for all $x,y \in \mS_2$
such that $x < y$ there exists $u \in \mS_2$ such that $u< y$ and $u \perp x$, and (b)
for any three incomparable elements 
of $\mS_2$ there is an element in $\mS_2$ that is larger than two out of the three, and incomparable to the third; 
\item \emph{nice} (following terminology from~\cite{DrosteHollandMacpherson}): for every $x,y \in \mS_2$ such that $x \perp y$ there exists $z \in \mS_2$ such that $z > x$ and $z \perp y$. 
\item \emph{without joins}:
for all $x,y,z \in \mS_2$ with $x,y \leq z$ and $x,y$ incomparable, there exists a $u \in \mS_2$
such that $x,y \leq u$ and $u < z$. 
\end{itemize}
It can be shown by a back-and-forth argument (and it also follows from results of Droste~\cite{Droste87Proc} and Droste, Holland, and Macpherson~\cite{DrosteHollandMacphersonII}) that all countable, dense, unbounded, nice, and binary branching semilinear orders without joins are isomorphic to $(\mS_2; \leq)$; see Proposition \ref{prop:backnforth} for details. 
 
Since all these properties of $(\mS_2;\leq)$ can be expressed by 
first-order sentences, it follows that $(\mS_2;\leq)$
is \emph{$\omega$-categorical}: it is, up to isomorphism, the unique countable model of its first-order theory.  
It also follows from general principles
that the first-order theory $T$ of $({\mathbb S}_2;\leq)$ is \emph{model complete}, that is, 
embeddings between models of $T$
 preserve all first-order formulas, and that $T$ is the 
 \emph{model companion} of the theory of semilinear orders, that is, has the same universal consequences; again, we refer to~\cite{Hodges} (Theorem 8.3.6). 


For $k \in \mathbb{N}$, a relational structure $\Delta$ is \emph{$k$ set--homogeneous} if whenever $A$ and $B$ are isomorphic $k$--element substructures of $\Delta$, there is an automorphism $g$ of $\Delta$ such that $g[A]=B$. In~\cite{Droste}, Droste studies 2 and 3 set--homogeneous semilinear orders. Of particular relevance here, Droste proved that $(\mS_2;\leq)$ is the unique countably infinite, non-linear, 3 set--homogeneous semilinear order (see Theorem 6.22 of~\cite{Droste}). 

The structure $(\mS_2;\leq)$ plays an important role in the study of a natural class of \emph{constraint satisfaction problems (CSPs)} in theoretical computer science. CSPs from this class have been studied in artificial intelligence for qualitative reasoning about branching time~\cite{Duentsch,Hirsch,BroxvallJonsson}, and, independently, in computational linguistics~\cite{Cornell,BodirskyKutz} under the name \emph{tree description} or \emph{dominance}  constraints. 

A \emph{reduct} of a relational structure $\Delta$ is a relational structure $\Gamma$ with the same domain as $\Delta$ such that every relation of $\Gamma$ has a first-order definition over $\Delta$ without parameters
(this slightly non-standard definition is common practice, see e.g.~\cite{RandomReducts,Thomas96,JunkerZiegler}). 
All reducts of a countable $\omega$-categorical structure
are again $\omega$-categorical~(Theorem 7.3.8 in~\cite{HodgesLong}). 
 In this article we study the reducts of $(\mS_2;\leq)$. Two structures $\Gamma$ and $\Gamma'$ with the same domain are called \emph{(first-order) interdefinable} when
$\Gamma$ is a reduct of $\Gamma'$, and $\Gamma'$ is a reduct of $\Gamma$. We show that the reducts $\Gamma$ of $(\mS_2;\leq)$ fall into three equivalence classes with respect to interdefinability: either
$\Gamma$ is interdefinable with $(\mS_2;=)$,
with $(\mS_2;\leq)$, or with 
$(\mS_2;B)$, where $B$ is the ternary \emph{betweenness relation}. The latter relation is defined by 
$$B(x,y,z) \; \Leftrightarrow \; (x < y < z) \vee (z < y < x) \vee (x < y \wedge y \perp z) \vee (z < y \wedge y \perp x) \; .$$

We also classify the \emph{model-complete cores} of the reducts of $(\mS_2;\leq)$. A structure $\Gamma$ is called \emph{model complete} \dff its
first-order theory is model complete. 
A structure $\Delta$ is a \emph{core} \dff all endomorphisms of $\Delta$ are embeddings. 
It is known that every $\omega$-categorical
structure $\Gamma$ is \emph{homomorphically equivalent}
to a model-complete core $\Delta$ (that is, there is a homomorphism from $\Gamma$ to $\Delta$ and vice versa; see~\cite{Cores-journal,BodHilsMartin}).
The structure $\Delta$ is unique up to isomorphism, $\omega$-categorical, and called
the \emph{model-complete core} of $\Gamma$. We show that for every reduct $\Gamma$ of $(\mS_2;\leq)$, 
the model-complete core 
of $\Gamma$ is
interdefinable with precisely one out of a list of ten structures (Corollary~\ref{cor:mc-cores}). 
The concept of model-complete cores is important for the aforementioned applications in constraint satisfaction, and implicitly used 
in complete complexity classifications for
the CSPs of reducts of $({\mathbb Q};<)$ 
and the CSPs of reducts of the random graph~\cite{tcsps-journal,BodPin-Schaefer-both};
also see~\cite{Bodirsky-HDR}. 
Our results have applications in this context which will be described in Section~\ref{sect:csp}.

There are alternative formulations of our results in the language of permutation groups and transformation monoids, which also plays an important role 
in the proofs.   
By the theorem of Ryll-Nardzewski (see, e.g., Corollary 7.3.3.~in Hodges~\cite{Hodges}), 
two $\omega$-categorical structures
are first-order interdefinable if and only if they have the
same automorphisms. 
Our result about the reducts of $(\mS_2;\leq)$
up to first-order interdefinability is equivalent
to the statement that there
are precisely three subgroups of 
$\Sym(\mathbb S_2)$
 that contain the automorphism group of 
$(\mS_2;\leq)$ and that are closed  with respect to the \emph{topology of pointwise convergence}, i.e., the product topology on $(\st)^{\mathbb S_2}$ where $\mathbb S_2$ is taken to be discrete.
The link to transformation monoids comes from
the fact that a countable $\omega$-categorical structure $\Gamma$ is model complete if and only if $\Aut(\Gamma)$ is dense in the monoid $\Emb(\Gamma)$ of self-embeddings of $\Gamma$, i.e., the closure $\overline{\Aut(\Gamma)}$ of $\Aut(\Gamma)$ in $(\st)^{\mathbb S_2}$ equals $\Emb(\Gamma)$; see~\cite{RandomMinOps}. Consequently, $\Gamma$ is a model-complete core if and only if $\Aut(\Gamma)$ is dense in the endomorphism monoid $\End(\Gamma)$ of $\Gamma$, i.e., $\overline{\Aut(\Gamma)}=\End(\Gamma)$.

The proof method for showing 
our results relies on
an analysis of the endomorphism monoids of reducts
of $(\mS_2;\leq)$. For that, we use a Ramsey-type
statement for semilattices, due to Leeb~\cite{Lee-vorlesungen-ueber-pascaltheorie} (cf.~also~\cite{GR-some-recent-developments}).
By results from~\cite{BP-reductsRamsey,BPT-decidability-of-definability}, that statement implies that if a reduct 
of $(\mS_2;\leq)$ has an endomorphism
that does not preserve a relation $R$, then
it also has an endomorphism that does not preserve $R$ and that behaves \emph{canonically} in a formal sense defined in Section~\ref{sect:prelims}. Canonicity
allows us to break the argument into finitely many cases.

We also mention a conjecture of Thomas,
which states that every countable homogeneous structure
$\Delta$ 
with a finite relational signature has only finitely many reducts up to interdefinability~\cite{RandomReducts}. 
By \emph{homogeneous} we mean here
 that every isomorphism between finite
substructures of $\Delta$ can be extended to an automorphism of $\Delta$. 
Thomas' conjecture has
been confirmed for various fundamental 
homogeneous structures, with particular activity in recent years~\cite{Cameron5,RandomReducts,Thomas96,Bennett-thesis,JunkerZiegler,Pon11,Poset-Reducts,42,LinmanPinsker,BodJonsPham}. 
The structure $(\mS_2;\leq)$ is not homogeneous, but interdefinable with a homogeneous structure with a finite relational signature, so it falls into the scope of Thomas' conjecture.

\section{Statement of Main Results}
\label{sect:results}
\begin{figure}
\begin{center}
\includegraphics[scale=.4]{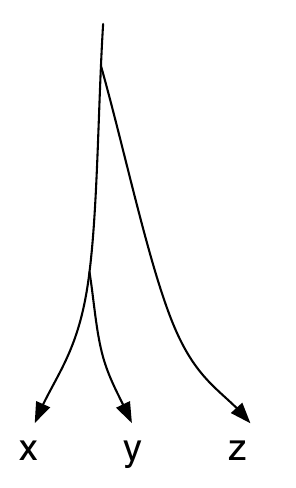} 
\end{center}
\caption{Illustration of $C(z,xy)$.}
\label{fig:c}
\end{figure}

To state our classification result, we need
to introduce some 
homogeneous structures that
appear in it. 
We have mentioned that
$(\mS_2;\leq)$ is not homogeneous, 
but interdefinable with a homogeneous structure
with finite relational signature. Indeed, to obtain
a homogeneous structure
we can add a single first-order definable ternary relation $C$ to $(\mS_2;\leq)$, defined as
\begin{align}
C(z,xy) \quad :\Leftrightarrow \quad x \perp y \; \wedge \;\exists u (x < u \wedge y < u \wedge u \perp z) \; . 
\label{C}
\end{align}
See Figure~\ref{fig:c}.

We omit the comma between the last two arguments of $C$ on purpose, since it increases readability, pointing out the symmetry $\forall x,y,z 
\; (C(z,xy) \Leftrightarrow C(z,yx))$. 
As pointed out in \cite{DrosteHollandMacphersonI}, it follows from Theorem 5.31 in~\cite{Droste} that the structure $(\mS;\leq,C)$ is homogeneous (for details see Proposition \ref{prop:backnforth}). 
Clearly, $(\mS_2;\leq)$ and $(\mS_2;\leq,C)$ are interdefinable. 

We write $(\mL_2;C)$ for the structure induced in $(\mS_2;C)$ by any maximal antichain of $(\mS_2;\leq)$.
It is straightforward
to verify that $(\mL_2;C)$ satisfies the axioms C1-C8 
given in~\cite{BodJonsPham}, and hence is isomorphic to
the \emph{homogeneous binary branching $C$-relation on leaves} which 
is also denoted by $(\mL_2;C)$ in~\cite{BodJonsPham} (see Lemma 3.8 in~\cite{BodJonsPham}). 
The reducts of $(\mL_2;C)$ were classified in~\cite{BodJonsPham}. 
We mention in passing that the structure $(\mL_2;C')$, where $C'(x,y,z) \Leftrightarrow \big(C(x,yz) \vee (y=z \wedge x \neq y)\big)$, is a so-called \emph{C-relation}; we refer to~\cite{AdelekeNeumann} for the definition since we will not make further use of it. 


It is known that two $\omega$-categorical structures have the same endomorphisms
if and only if they are existentially positively
interdefinable, that is, if and only if each relation in one of the structures can be defined by an existential positive formula in the other structure~\cite{RandomMinOps}. We can now state one of our main results. 

\begin{theorem}\label{thm:4case}
Let $\Gamma$ be a reduct of $(\mS_2; \leq)$. Then at least one of the following cases applies.
\begin{itemize}
\item[(1)] $\End(\Gamma)$ contains a function whose range induces a chain in $(\st;\leq)$, and 
$\Gamma$ is homomorphically equivalent to a reduct of the order of the rationals $(\mQ; <)$.
\item[(2)] $\End(\Gamma)$ contains a function whose range induces an antichain in $(\st;\leq)$, and $\Gamma$ is homomorphically equivalent to a reduct of $(\mL_2; C)$.
\item[(3)] 
$\End(\Gamma)$ equals $\overline{\Aut(\mathbb S_2; B)}$; equivalently,
$\Gamma$ is existentially positively interdefinable
with $(\mS_2;B)$. 
\item[(4)] $\End(\Gamma)$ equals $\overline{\Aut(\mathbb S_2; \leq)}$; equivalently,
$\Gamma$ is existentially positively interdefinable
with $(\mS_2;<,\perp)$. 
\end{itemize}
\end{theorem}

The reducts of
$(\mL_2;C)$ have been classified in~\cite{BodJonsPham}. Each reduct of $(\mL_2;C)$ is interdefinable
with either
\begin{itemize}
\item $(\mL_2;C)$ itself,
\item $(\mL_2;D)$ where $D(x,y,u,v)$ has the first-order definition $$(C(u,xy) \wedge C(v,xy)) \vee (C(x,uv) \wedge C(y,uv))$$ over $(\mL_2;C)$, or
\item $(\mL_2;=)$.
\end{itemize}


The reducts of $(\mQ;<)$ have been classified in~\cite{Cameron5}. In order
to keep the formulas compact, we write 
$\overrightarrow{x_1\cdots x_n}$ 
whenever $x_1,\ldots,x_n\in \mathbb{Q}$ are such that $x_1<\cdots<x_n$. 
Cameron's theorem states that each reduct of $(\mQ;<)$ is interdefinable with either
\begin{itemize}
\item the dense linear order $(\mathbb{Q};<)$ itself,
\item the structure $(\mathbb{Q}; \Betw)$, where $\Betw$ is the ternary relation 
$$\big \{(x,y,z) \in {\mathbb Q}^3 : \overrightarrow{xyz} \, \vee \,
\overrightarrow{zyx} \big \} \, ,$$
\item the structure $(\mathbb{Q}; \Cyc)$, where $\Cyc$ is the ternary
  relation 
   $$\big \{(x,y,z) : \overrightarrow{xyz} \vee \overrightarrow{yzx} \vee \overrightarrow{zxy} \big \} \, ,$$
\item the structure $(\mathbb{Q}; \Sep)$, where $\Sep$ is the 4-ary relation 
\begin{align*}
 \big \{(x_1,y_1,x_2,y_2) : \; &
 \overrightarrow{x_1x_2y_1y_2} \vee \overrightarrow{x_1y_2y_1x_2} \vee
  \overrightarrow{y_1x_2x_1y_2} \vee \overrightarrow{y_1y_2x_1x_2} \\
  \vee \; & \overrightarrow{x_2x_1y_2y_1}  \vee
  \overrightarrow{x_2y_1y_2x_1} \vee \overrightarrow{y_2x_1x_2y_1} \vee \overrightarrow{y_2y_1x_2x_1} \big \} \, , \text{ or }
  \end{align*}
\item the structure $(\mathbb{Q};=)$.
\end{itemize}

\begin{corollary}\label{cor:mc-cores}
Let $\Gamma$ be a reduct of $(\mS_2; \leq)$. Then its model-complete core has only one element, or is 
interdefinable with $(\mS_2;<,\perp)$,
 $(\mS_2;B)$, $(\mL_2;C)$, $(\mL_2;D)$, 
 $(\mQ;<)$, $(\mQ;\Betw)$, $(\mQ;\Cyc)$, 
 $(\mQ;\Sep)$, or $(\mQ;\neq)$. 
\end{corollary}



\begin{theorem}\label{thm:groups}
Let $\Gamma$ be a reduct of $(\mS_2; \leq)$. 
Then $\Gamma$ is first-order
interdefinable with $(\mS_2; \leq)$, $(\mS_2; B)$, or $(\mS_2; =)$. Equivalently, $\Aut(\Gamma)$ equals either 
$\Aut(\mS_2; \leq)$, $\Aut(\mS_2; B)$, or $\Aut(\mS_2; =)$.
\end{theorem}

The closed subgroups of 
$\Sym(\mathbb S_2)$ are precisely the automorphism groups
of structures with domain $\mathbb S_2$ (see, e.g., \cite{Oligo}). 
Moreover, the closed subgroups of $\Sym(\mathbb S_2)$ that contain $\Aut(\mS_2;\leq)$
are precisely the automorphism groups 
of reducts of $(\mS_2;\leq)$, becuase $(\mS_2;\leq)$ is $\omega$-categorical; again, see~\cite{Oligo} for background.  Therefore, the following is an immediate consequence of Theorem~\ref{thm:groups}.

\begin{corollary}
The closed subgroups of $\Sym(\mathbb S_2)$ containing $\Aut(\mS_2; \leq)$ are precisely the permutation groups $\Aut(\mS_2;\leq)$, $\Aut(\mS_2; B)$, and $\Aut(\mS_2; =)$. 
\end{corollary}

\section{Preliminaries}
\label{sect:prelims}
In the introduction we gave an explicit first-order axiomatisation of $(\mS;\leq)$. Although this follows from results in \cite{Droste87Proc} and \cite{DrosteHollandMacphersonII}, we provide details here for the convenience of the reader; also proving  the claim about the homogeneity of $(\mS_2;\leq,C)$ made in Section~\ref{sect:results}. 
We then review the Ramsey properties of $(\mS_2;\leq)$ after the expansion with a suitable linear order in Section~\ref{sect:ramsey}. The Ramsey property will be used in our proof via the concept of \emph{canonical functions}; they will be introduced in 
 Section~\ref{sect:canonical}. 

\subsection{Homogeneity of 
$({\mathbb S}_2;\leq,C)$}
\label{sect:s2}
We show that all countable semilinear orders that are dense, unbounded, binary branching, nice and without joins are isomorphic. Since drafting this paper, we have learnt that this follows as a special case of Proposition 2.7 of~\cite{DrosteHollandMacphersonII} and Theorem 4.4 of~\cite{Droste87Proc}. The axioms we have given explicitly correspond to their notion of an almost normal tree of type $(1,(0,0),\{2\})$. 
For completeness, we provide a self-contained proof which also establishes the homogeneity of $(\mS;\leq,C)$; though as stated in the introduction of \cite{DrosteHollandMacphersonI}, 
this follows from Theorem 5.31 
of~\cite{Droste}.


For subsets $U,V$ of a poset, 
we write $U < V$ if 
$u<v$ holds for all $u \in U$ and $v \in V$.
The notation $U \leq V$ and $U \perp V$ is defined
analogously. We also write $u < V$ for $\{u\} < V$
and $u \perp V$ for $\{u\} \perp V$. 

\begin{lemma}
\label{lem:key}
 Let $(P;\leq)$ be a dense, unbounded, nice, and 
 binary branching semilinear order without joins. 
 Let $U,V,W \subseteq P$ be finite subsets such that $U$ is non-empty, 
 $U < V$, $U \perp W$, and $C(w,u_1u_2)$ for all $w \in W$ and incomparable $u_1,u_2 \in U$. Then there exists an 
 $x \in P$ such that  $U < x$, $x<V$, and $x \perp W$. 
\end{lemma}
\begin{proof}
First note that if $V \cup W$ is empty, the lemma follows from (upward) unboundedness of $P$. So assume throughout that $V\cup W$ is non-empty.
For $p,q \in V \cup W$, 
define $p \lhd q$ if 
 \begin{itemize}
 \item $p < q$,
 \item $p \perp q$ and $u < p$ for all $u \in U$, or
 \item $C(q,pu)$ for all $u \in U$. 
\end{itemize}
Note that $\lhd$ is a strict partial order on $V \cup W$: irreflexivity is immediate from the definition; 
to verify
transitivity, let $a,b,c \in V \cup W$
such that $a \lhd b$ and $b \lhd c$ and check the various configurations of $a,b,c$ with respect to the (non-empty) set $U$.
First assume that $c \perp b$. 
If $b > U$, then $b \in V$ and so $a \lhd b$ implies that $a < b$, 
thus $c \perp a$ from the semilinearity of $<$. Moreover, either $u < a$ for all $u \in U$, or $C(c,au)$ is witnessed by $b$ for all $u \in U$; hence 
$a \lhd c$.
Otherwise, suppose that $C(c,bu)$ for all $u \in U$. 
If $a<b$ then also $C(c,au)$ for all $u \in U$, implying $a \lhd c$. 
If $a \perp b$ and $u < a$ for all $u \in U$, then also $a \perp c$; so $a \lhd c$.
If $C(b,au)$ for all $u \in U$, then also $C(c,au)$ for all $u \in U$, also yielding $a \lhd c$.
Assume now that we have $b < c$.
If $a<b$, then $a<c$, by the transitivity of $<$, so $a \lhd c$. Moreover, if $b \in V$ then $a \lhd b$ implies $a<b$, so $a<b<c$ and we are done.
So suppose that $a \perp b$ and $b \in W$.
If $c \in V$ then either $a \perp b$ and $U < a$, or $C(b,au)$ for all $u \in U$; in either case we have $a<c$, hence $a \lhd c$.
Instead if $c \in W$ then we either have that $C(b,au)$ for all $u \in U$, in which case $C(c,au)$ for all $u \in U$, or we have that $u < a$ for all $u \in U$, in which case $a \perp c$; in either case, $a \lhd c$.

Therefore, as $V \cup W$ is non-empty, there exists an element $m \in V \cup W$ that is minimal with respect to $\lhd$. 

We prove the statement of the lemma by induction on the
number of elements of $U$. 
Since $U$ is non-empty and finite, it contains a maximal element $u_0$ with respect to $<$.
If there is just one such element, we distinguish whether
$m \in V$ or $m \in W$. If $m \in V$ then we choose
$x \in P$ such that $u_0 < x < m$; such an $x$ exists by density of $(P;\leq)$.
The minimality of $m$ with respect to $\lhd$ implies that $m \perp W$ and $m$ is the minimum of $V$, 
as $V$ is linearly ordered by $\lhd$. 
 So by transitivity of $<$, we have $x \perp W$ and $U < x < V$, as $u_0$ is the only maximal element in $U$.
If $m \in W$ then we choose $x \in P$ such that 
$u_0 < x$
and $x \perp m$; such an $x$ exists since $(P;\leq)$ is nice. As before, we
 have $U < x$. Moreover, $x \perp W$
 and $x < V$ hold by the minimality of $m$ with respect to $\lhd$. 

Now consider the case that there are two distinct maximal elements $u_0,u_1 \in U$. 
Again we distinguish two subcases. If $m \in V$
then there exists an element $x \in P$ such that
$u_0,u_1 < x$ and $x<m$, 
since $(P;\leq)$ is without joins. 
We then have that $U < x$ since 
$u_0,u_1 < x$ are the only two maximal elements of $U$. Moreover, 
$x < V$ since $x < m$, and $x \perp W$ 
by minimality of $m$ with respect to $\lhd$. Otherwise, $m \in W$. 
Since we have 
$C(m,u_0u_1)$ by assumption, there exists an element
$x \in P$ such that $x > u_0,u_1$ and $x \perp m$,
and this element $x$ satisfies the required conditions: 
 $x < V$ and $x \perp W$ by the minimality of $m$ with respect to $\lhd$ and clearly $x>U$.

Now suppose that there are at least three distinct maximal
elements $u_0,u_1,u_2$ in $U$. 
Since $(P;\leq)$ is binary branching, 
there is
an $s \in P$ larger than two out of $u_0,u_1,u_2$ and incomparable to the third; without loss of generality say that $s > u_0, u_1$ and $s \perp u_2$. 
Note that $s < V$ and that
$C(w,us)$ for every $w \in W$ and 
every $u \in U$ incomparable with $s$, which implies that $s$ is incomparable with $W$.
Hence, we can apply the inductive assumption for 
the non-empty 
set $U' := U \cup \{s\} \setminus \{u_0,u_1\}$ instead of $U$,
which has one element less than $U$. The element $x \in P$ that we obtain for $U'$ also satisfies the requirements that we have for $U$: we have $x < V$ and $x \perp W$, and $x > U$ follows
 from $x > U'$ since $x > s$ implies $x > u_0,u_1$. 
\end{proof}

\begin{proposition}\label{prop:backnforth}
All countable semilinear orders that are
dense, unbounded, binary branching, nice, and without joins are isomorphic to $({\mathbb S}_2;\leq)$. The structure $({\mathbb S}_2;\leq,C)$ is homogeneous. 
\end{proposition}

\begin{proof}
The proof uses a standard back-and-forth argument, where we inductively construct an isomorphism between 
two semilinear orders $(P;\leq)$ and $(Q;\leq)$ that satisfy the properties given in the statement, by alternating between steps that make sure that the function will be defined everywhere (going forth) and steps that make sure that the function will be a surjection (going back). 
Let $\Gamma$ and
$\Delta$ be the expansions of 
$(P;\leq)$ and $(Q;\leq)$
 with the signature $\{\leq,C\}$ where $C$ denotes the relation as defined in (\ref{C}) at the beginning of Section~\ref{sect:results}. 

We fix enumerations $(p_i)_{i \in \omega}$ and $(q_j)_{j \in \omega}$ of $P$ and $Q$, respectively. Assume that $D \subseteq P$ is a finite subset of $P$ and that $\rho \colon D \rightarrow E$ is an isomorphism between
the substructure induced by $D$ in $\Gamma$ and the substructure induced by $E$ in $\Delta$. 
Let $k \in \omega$ be smallest such that $p_k \in P \setminus D$. To go forth we need to extend the domain of the partial isomorphism $\rho$ to $D \cup \{p_k\}$. 
Let $D_> := \{a \in D : a > p_k\}$ and $D_< := \{a \in D : a < p_k\}$ and $D_\perp := \{a \in D : a \perp p_k\}$. 
In each case we describe the element $q \in Q$ 
such that $\rho(p_k) := q$ defines an extension 
of $\rho$ which is a partial isomorphism between
$(P;\leq,C)$ and $(Q;\leq,C)$. 




\paragraph{\textbf{Case 1: $D_<$ is empty.}}
If $D_> \cup D_{\perp}$ is also empty, any $q \in Q$ will suffice for the image of $p_k$. So we assume that $D_> \cup D_{\perp}$ is non-empty.
Suppose first that there is an element
$v \in D_>$ such that $v \perp w$ for all $w \in D_\perp$; choose $v$ minimal
with these properties. 
In this case we can choose $q \in Q$ such that 
$q < \rho(v)$ by the unboundedness of $(Q;\leq)$. Then $q \perp \rho[D_\perp]$ by the transitivity of $<$ and $q < \rho[D_>]$ by the minimality of $\rho(v)$ 
in $\rho[D_>]$. 
Moreover, it is clear from the definition of $C$ that for any $w_1, w_2 \in D_{\perp}$ we have $C(p_k,w_1 w_2) \Leftrightarrow C(v,w_1 w_2) \Leftrightarrow C(\rho(v),\rho(w_1) \rho(w_2)) \Leftrightarrow C(q,\rho(w_1) \rho(w_2))$ and $C(w_1, p_k w_2) \Leftrightarrow C(w_1, v w_2) \Leftrightarrow C(\rho(w_1), \rho(v) \rho(w_2)) \Leftrightarrow C(\rho(w_1), q \rho(w_2))$, so the extension of $\rho$ indeed 
yields a partial isomorphism.

Otherwise, there exists an element $w_0 \in D_\perp$
such that $w_0 < v$ for all $v \in D_>$. 
Choose $w_0$ minimal with respect
to the relation $\lhd$ as defined in the proof of Lemma~\ref{lem:key} for $V := D_>$, 
$W := D_\perp$, and $U := \{p_k\}$
(clearly, we then have $U<V$ and $U \perp W$ while the condition in Lemma~\ref{lem:key} that $C(w,u_1u_2)$ for all $w \in W$ and incomparable $u_1,u_2 \in U$ becomes void since $|U| = 1$).  
Now partition $D_{\perp}$ into $\overline{U}:=\{u \in D_{\perp} \mid C(p_k, u w_0) \textrm{ or } u \geq w_0\}$ and $\overline{W} := D_{\perp} \setminus \overline{U}$.
By Lemma 3.1 applied to $U := \rho[\overline{U}], V:= \rho[D_{>}],$ and $W := \rho[\overline{W}]$ (which satisfy the assumptions of Lemma 3.1: clearly $U<V$, and $U \perp W$ follows from $C(w, u p_k)$ for all $w\in \overline{W}$ and $u \in \overline{U}$,  
while for any incomparable $\rho(u_1),\rho(u_2) \in U$ we have $C(p_k,u_1 u_2)$, 
so $C(w,u_1 u_2)$ for any $\rho(w) \in W$ by the definition of $W$, then $C(\rho(w),\rho(u_1) \rho(u_2))$ since $\rho$ is a partial isomorphism),
we obtain 
an element $x \in Q$ such that 
$U < x$, $x < V$, and $x \perp W$. 
Another application of this lemma,
 this time applied to $V := \{x\} \cup \rho[D_>]$ and $U$ and $W$ as before
gives us an element
$x' \in Q$ with the same properties and $x' < x$. 
Since $(Q;\leq)$
is binary branching there exists an element
$q \in Q$ with $q < x$ and $q \perp x'$. 
To see that $q$ has the required properties so that the extension of $\rho$ is an isomorphism, first note that $q \perp \rho[D_{\perp}]$ and $q<\rho[D_>]$. Furthermore, for any $\rho(u_1), \rho(u_2) \in U$, we have $C(q,\rho(u_1) \rho(u_2))$ witnessed by $x'$, while for any $\rho(u) \in U$ and $\rho(w) \in W$ we have $C(\rho(w),\rho(u) q)$ witnessed by $x$. Finally, note that for any $\rho(w_1), \rho(w_2) \in W$ we have 
\begin{align*}
C(q,\rho(w_1) \rho(w_2)) 
& \Leftrightarrow C(x,\rho(w_1) \rho(w_2)) 
\Leftrightarrow C(\rho(w_0),\rho(w_1) \rho(w_2))
\end{align*} 
and 
\begin{align*}
C(\rho(w_2),\rho(w_1) q) \Leftrightarrow C(\rho(w_2), \rho(w_1) x) \Leftrightarrow C(\rho(w_2), \rho(w_1) \rho(w_0))
\end{align*} so, as 
$C(w_0, w_1 w_2) \Leftrightarrow C(p_k,w_1, w_2)$,  $C(w_2, w_1 w_0) \Leftrightarrow C(w_2, w_1 p_k)$, and
$\rho$ is assumed to be a partial isomorphism on $D$, 
we have 
$C(p_k,w_1 w_2) \Leftrightarrow C(q, \rho(w_1) \rho(w_2))$ 
and $C(w_2,w_1 p_k) \Leftrightarrow C(\rho(w_2),\rho(w_1) q)$.

\paragraph{\textbf{Case 2: $D_<$ is non-empty.}}
We apply Lemma~\ref{lem:key} to
$U:= \rho[D_<]$, $V := \rho[D_>]$, and $W := \rho[D_\perp]$. 
The element $x$ from the statement of 
Lemma~\ref{lem:key} has the properties that we require for $q$, namely $q \perp \rho[D_\perp]$,
$q < \rho[D_>]$, and $q > \rho[D_<]$. 
Moreover, for any $\rho(w_1), \rho(w_2) \in W$ and $u \in U$ we have 
\begin{align*}
C(q, \rho(w_1)\rho(w_2)) & \Leftrightarrow C(\rho(u), \rho(w_1)\rho(w_2)) 
\Leftrightarrow C(u,w_1 w_2) \Leftrightarrow C(p_k,w_1 w_2)
\end{align*}
 and 
 \begin{align*}
 C(\rho(w_2),\rho(w_1) q) \Leftrightarrow C(\rho(w_2),\rho(w_1) \rho(u)) \Leftrightarrow C(w_2,w_1 u) \Leftrightarrow C(w_2,w_1 p_k)
 \end{align*}
so the extension of $\rho$ yields a partial isomorphism. 

This allows us to take the step going forth. 
To take the step going back, we need to extend the range of $\rho$ to $D' \cup \{q_k\}$ where $k$ is the first such that $q_k \in Q \setminus D'$. The argument is analogous to the argument given above for going forth. 
This concludes the back-and-forth and the result follows.
\end{proof}

\subsection{The convex linear Ramsey extension}
\label{sect:ramsey}

Let $(S;\leq)$ be a semilinear order. A linear order $\prec$ on $S$ is called a \emph{convex linear extension of $\leq$} \dff the following two conditions hold; here, the relations $<$ and $C$ are defined over $(S;\leq)$ as they were defined over $(\mS_2;\leq)$.
\begin{itemize}
\item $\prec$ is an extension of $<$, i.e., $x<y$ implies $x\prec y$ for all $x,y\in S$;
\item for all  $x,y,z\in S$ we have that $C(x,yz)$ implies that $x$ cannot lie between $y$ and $z$ with respect to $\prec$, i.e., 
$(x \prec y\wedge x\prec z)\vee (y\prec x\wedge z\prec x)$.
\end{itemize}
For finite semilinear orders $(S;\leq)$, the convex linear extensions are precisely the 
linear orders $\prec$ that can be defined recursively as follows. There exists a largest element $r \in S$; 
let $v_1,\dots,v_s$ be the maximal elements below $r$. For each $i \leq s$, we define $\prec$ recursively on the semilinear order induced by $S_i := \{u \in S \mid u < v_i\}$ in $(S;\leq)$. Note that $\{r\},S_1,\dots,S_s$ partition 
$S$, and we finally put $S_1 \prec \cdots \prec S_s \prec \{r\}$. 

Using Fra\"{i}ss\'{e}'s theorem~\cite{HodgesLong} one can show that in the case of $(\mS_2;\leq)$, there exists a convex linear extension $\prec$ of $\leq$ such that $(\mS_2;\leq,C,\prec)$ is homogeneous and such that $({\mathbb S}_2;\leq,\prec)$ is \emph{universal} in the sense that it contains all isomorphism types of convex linear extensions of finite semilinear orders; this extension is unique in the sense that all expansions of $(\mS_2;\leq,C)$ by a convex linear extension with the above properties are isomorphic. We henceforth fix any such extension $\prec$. 
The structure $(\mS_2;\leq,C,\prec)$
is homogeneous and therefore also model complete. Moreover, the structure is combinatorially well-behaved in the following sense. For  structures $\Sigma, \Pi$ in the same language, we write $\mix \Sigma \Pi$ for the set of all 
embeddings of $\Pi$ into $\Sigma$.

\begin{defn}\label{def:ramseystructure}
A countable homogeneous relational structure $\Delta$ is called a \emph{Ramsey structure} \dff for all finite substructures $\Omega$ of $\Delta$, all substructures $\Gamma$ of $\Omega$, and all $\chi\colon\mix \Delta \Gamma\To 2$ there exists an $e_1 \in \mix \Delta \Omega$ such that $\chi$ is
constant on $e_1 \circ \mix {\Omega} \Gamma$ (which denotes the set of compositions of $e_1$ with a function from $\mix \Omega \Gamma$). A countable $\omega$-categorical structure $\Delta$ 
is called \emph{Ramsey} if 
the (necessarily homogeneous) 
relational structure
whose relations are precisely the first-order definable relations in $\Delta$ is Ramsey. 
\end{defn}

The following theorem is a special case of a Ramsey-type statement for semilinearly ordered semilattices due to Leeb~\cite{Lee-vorlesungen-ueber-pascaltheorie} (also see~\cite{GR-some-recent-developments}, page 276). A \emph{semilinearly ordered semilattice} $(S;\vee,\leq)$ is a semilinear order $(S;\leq)$ which is closed under the binary function $\vee$, the \emph{join} function, satisfying for all $x$ and $y$, that $x \vee y$ is the least upper bound of $\{x,y\}$ with respect to $\leq$. If $\prec$ is a convex linear extension of $\leq$, then $(S;\vee,\leq,\prec)$ is a convex linear extension of the semilinearly ordered semilattice  $(S;\vee,\leq)$. By Fra\"{i}ss\'{e}'s Theorem~\cite{HodgesLong} 
there is a countably infinite homogeneous structure $(\mathbb{T};\vee,\leq,\prec)$ 
which is the Fra\"{i}ss\'{e} limit of the class of finite, semilinearly ordered semilattices with a convex linear extension, as this class
is an amalgamation class. 

\begin{theorem}[Leeb]
\label{thm:Leeb}
 $(\mathbb{T};\vee,\leq,\prec)$ is a Ramsey structure.
\end{theorem}

\begin{corollary}
 $(\mS_2;\leq,C,\prec)$ is a Ramsey structure.
\end{corollary}

\begin{proof}
The same relations are first-order definable in $(\mathbb{T};\vee,\leq,\prec)$ and in $(\mathbb{T};\leq,\prec)$, and so Theorem \ref{thm:Leeb} above implies that 
$(\mathbb{T};\leq,\prec)$ is a Ramsey structure. Every finite substructure of $(\mS_2;\leq,\prec)$ is isomorphic to a substructure of $(\mathbb{T};\leq,\prec)$ and vice versa, so they have the same age and the two structures satisfy the same universal sentences. Hence, as $(\mS_2;\leq,\prec)$ is model complete, it is the model companion of $(\mathbb{T};\leq,\prec)$. 
Theorem~3.15 of \cite{BodirskyRamsey} states that the model companion of an $\omega$-categorical Ramsey structure is Ramsey, so we conclude that $(\mS_2;\leq,\prec)$ is a Ramsey structure, and so is the homogeneous structure $(\mS_2;\leq,C,\prec)$ because $C$ is first-order definable in $(\mS_2;\leq,\prec)$. 
\end{proof}

\subsection{Canonical functions.}
\label{sect:canonical}
The fact that $(\mS_2;\leq,C,\prec)$ is a relational homogeneous Ramsey structure implies that endomorphism monoids of reducts of this structure, and hence also of $(\mS_2;\leq,C)$, can be distinguished by  so-called \emph{canonical functions}.

\begin{defn} 
Let $\Delta$ be a structure, and let $a$ be an $n$-tuple of elements in $\Delta$. The \emph{type} of $a$ in $\Delta$ is the set of first-order formulas with free variables $x_1,\ldots,x_n$ that hold for $a$ in $\Delta$.
\end{defn}

\begin{defn} 
Let $\Delta$ and $\Gamma$ be structures. A \emph{type condition} between $\Delta$ and $\Gamma$ is a pair $(t,s)$, such that $t$ is the type on an $n$-tuple in $\Delta$ and $s$ is the type of an $n$-tuple in $\Gamma$, for some $n\geq 1$. A function $f\colon\Delta\to\Gamma$ \emph{satisfies} a type condition $(t,s)$ \dff the type of $(f(a_1),\ldots,f(a_n))$ in $\Gamma$ equals $s$ for all $n$-tuples $(a_1,\ldots,a_n)$ in $\Delta$ of type $t$. 

A \emph{behaviour} is a set of type conditions between $\Delta$ and $\Gamma$. We say that a function $f\colon\Delta\to\Gamma$ has a given behaviour \dff it satisfies all of its type conditions.
\end{defn}

\begin{defn} 
Let $\Delta$ and $\Gamma$ be structures. A function $f \colon \Delta\to\Gamma$ is \emph{canonical} \dff for every type $t$ of an $n$-tuple in $\Delta$ there is a type $s$ of an $n$-tuple in $\Gamma$ such that $f$ satisfies the type condition $(t,s)$. That is, canonical functions send $n$-tuples of the same type to $n$-tuples of the same type, for all $n\geq 1$.
\end{defn}

Note that any canonical function induces a function from the types over $\Delta$ to the types over $\Gamma$.

\begin{defn} Let $\mathcal{F} \subseteq (\mS_2)^{\mS_2}$. We say that $\mathcal{F}$ \emph{generates} a function $g\colon \mS_2\to \mS_2$ \dff $g$ is contained in the smallest closed (with respect to the topology of pointwise convergence) submonoid of $(\mS_2)^{\mS_2}$ which contains $\mathcal F$. The definition extends naturally to sets of functions being generated. 
\end{defn}

Note that $\mathcal{F}$ generates $g$
 \iff
for every finite subset $A\subseteq \mS_2$ there exists an $n\geq 1$ and $f_1,\ldots, f_n\in\mathcal{F}$ such that $f_1\circ\cdots\circ f_n$ agrees with $g$ on $A$ (see, e.g., Proposition~3.3.6 in~\cite{Bodirsky-HDR}).
 
Our proof relies on the following proposition which is a consequence of~\cite{BP-reductsRamsey, BPT-decidability-of-definability} and the fact that $(\mS_2;\leq,C,\prec)$ is a homogeneous Ramsey structure. For a structure $\Delta$ and elements $c_1,\ldots,c_n$ in that structure, let $(\Delta,c_1,\ldots,c_n)$ denote the structure obtained from $\Delta$ by adding the constants $c_1,\ldots,c_n$ to the language.

\begin{prop}\label{prop:canfcts}
Let $f\colon \mS_2\To \mS_2$ be any injective function, and let $c_1,\ldots,c_n\in \mS_2$. Then $\{f\}\cup\Aut(\mS_2;\leq,\prec)$ generates an injective function $g\colon \mS_2 \To \mS_2$ such that
\begin{itemize}
\item $g$ agrees with $f$ on $\{c_1,\ldots,c_n\}$;
\item $g$ is canonical as a function from $(\mS_2;\leq,C,\prec,c_1,\ldots,c_n)$ to $(\mS_2;\leq,C,\prec)$.
\end{itemize}
\end{prop}
\begin{proof}
Lemma 14 in~\cite{BPT-decidability-of-definability} proves the statement for all ordered homogeneous Ramsey structures with finite relational signature; $\mS_2$ is such a structure.
\end{proof}

\section{The Proof}
We start this section with a description of the functions in $\Emb(\mS_2;B)$ 
since they play an important role in the proof. Section~\ref{sect:ramsey-analysis}
contains the core of the classification
which is based on Ramsey theory. 
Our main result about endomorphism monoids of reducts of $(\mS_2;<)$, Theorem~\ref{thm:4case}, is shown in Section~\ref{sect:4case}. 
The classification of the
automorphism groups of
reducts of $(\mS_2;<)$, Theorem~\ref{thm:groups}, is not an immediate consequence of this result about the endomorphism monoids,
and we prove it in Section~\ref{sect:groups}.


\subsection{Rerootings and betweenness}
We start by examining what the automorphisms, self-embeddings, and endomorphisms of $(\mS_2;B)$ look like.

\begin{lemma}\label{lem:not-b}
Any function in $(\st)^{\st}$ that preserves $B$ is injective and preserves $\lnot B$. 
\end{lemma}
\begin{proof}
The existential positive formula 
$$(a=b)\vee (b=c)\vee (c=a)\vee \exists x (B(a,x,b)\wedge B(b,x,c))
$$ is equivalent to $\lnot B(a,b,c)$. 
Moreover, for all $a,b\in \st$ we have that $a\neq b$ \iff there exists $c\in\st$ such that $B(a,b,c)$, so inequality has an existential positive definition from $B$, and functions preserving $B$ must be injective. Hence, every endomorphism of $(\mS_2; B)$ is an embedding (cf.~the discussion in the introduction).
\end{proof}

\begin{defn}\label{defn:rerooting}
A \emph{rerooting} of $(\mS_2;<)$ is an injective function $f\colon\st\To\st$ for which there exists a set $S\subseteq \st$ such that
\begin{itemize}
\item $S$ is an upward closed chain, i.e., if $x\in S$ and $y\in\st$ satisfy $y>x$, then $y\in S$;
\item $f$ reverses the order $<$ on $S$;
\item $f$ preserves $<$ and $\perp$ on $\st\setminus S$;
\item whenever $x\in \st\setminus S$ and $y\in S$, then $x<y$ implies $f(x)\perp f(y)$ and $x\perp y$ implies $f(x)<f(y)$.
\end{itemize}
We then say that $f$ is a \emph{rerooting with respect to $S$}.
\end{defn}

It is not hard to see that whenever $S\subseteq \st$ is as above, then there is a rerooting with respect to $S$: it suffices to verify that the relation $<$
on the image of $f$ given by the conditions above is a partial order 
and that there are no elements $a,b,c \in \mS_2$ such that $f(a) < f(b)$, $f(a)<f(c)$, and $f(b) \perp f(c)$ (which would violate semilinearity). 
A rerooting with respect to $S$ is a self-embedding of $(\mS_2;<)$ if and only if $S$ is empty. 

The image of any rerooting with respect to $S$ is isomorphic to $(\mS_2;<)$ if and only if $S$ is a maximal chain or empty:
if $S$ is a chain that is not maximal and $f$ a rerooting, then there is some $a < S$. Then  
$f(a) \perp f(S)$ and hence
$\{f(a)\} \cup f(S)$ has no upper bound in the image of $f$; the image is not a semilinear order. Whereas it can be verified that the image of a rerooting with respect to a maximal chain $S$ is a dense, unbounded, binary branching, nice, semilinear order without joins (for each property one can pull back any instance of the universally bound quantifier via the inverse of the rerooting and the existence of the required element in the image is asserted by one of these properties in the pre-image) so by Proposition \ref{prop:backnforth} the image is isomorphic to $(\mS_2;<)$.
In particular, there exist rerootings which are permutations of $\st$ and which are not self-embeddings of $(\mS_2;<)$.

\begin{proposition}\label{prop:rerooting}
$\Emb(\mS_2;B)$ consists precisely of the rerootings of $(\mS_2;<)$.
\end{proposition}
\begin{proof}
To see that any rerooting preserves $B$, take $(a,b,c) \in B$ and $f$ a rerooting with respect to $S$ as in Definition \ref{defn:rerooting}. Without loss of generality, either $a<b<c$ or $c \perp a<b \perp c$. If $a < b <c$ then, depending precisely on the cardinality of $S \cap \{a,b,c\}$, either $f(a)<f(b)<f(c)$, or $f(c)\perp f(a) <f(b) \perp f(c)$, or $f(a) \perp f(c) < f(b) \perp f(a)$, or $f(c) < f(b) < f(a)$. Otherwise, assume that $c \perp a< b \perp c$. If $S$ omits only one element, it omits $c$, so $f(c)<f(b)<f(a)$. If $S$ omits two, they are $b$ and $a$, or $a$ and $c$; in which case $f(a)<f(b)<f(c)$, or $f(a)\perp f(c)<f(b) \perp f(a)$, respectively. If $S$ omits all three, then it behaves like the identity, so $f(c) \perp f(a) < f(b) \perp f(c)$. 
Then Lemma~\ref{lem:not-b} implies that 
$f \in \Emb(\mS_2;B)$. 

Conversely, suppose that
$f\in \Emb(\mS_2;B)$. We claim that $f\in\Emb(\st;<)$ or there exist $x,y\in\st$ such that $x<y$ and $f(x)>f(y)$. Suppose that
$f \notin \Emb(\mS_2;<)$. Then
 $f$ violates $\perp$ or $f$ violates $<$. 
Suppose $f$ violates $\perp$. Pick $a,b\in\st$ with $a\perp b$ and such that $f(a)<f(b)$. There exists $c\in\st$ such that $c>b$ and such that $B(a,c,b)$. Since $f$ preserves $B$ we then must have $f({c})<f(b)$, and our claim follows. Now suppose $f$ violates $<$, and pick $a,b\in\st$ with $a<b$ witnessing this. Then for any $c\in\st$ with $c>b$ we have $f({c})<f(b)$, as $f$ preserves $B$, proving the claim.

Let $S:=\{x\in\st\mid \exists y\in\st (x<y\wedge f(y)<f(x))\}$. By the above, we may assume that $S$ is non-empty. Since $f$ preserves $B$, it follows easily that whenever $x\in S$, $y\in \st$ and $x<y$, then $f(y)<f(x)$. From this and again because $f$ preserves $B$ it follows that $S$ is upward closed. Hence, $S$ cannot contain incomparable elements $x,y$, as otherwise for any $z\in S$ with $x<z$ and $y<z$ we would have $f(x)>f(z)$ and $f(y)>f(z)$, and so $f(x)$ and $f(y)$ would have to be comparable. But then $f$ would violate $\neg B$ on $\{x,y,z\}$.
So this $S$ satisfies the first part of Definition~\ref{defn:rerooting} and $f$ behaves on $S$ as required by the second part of the definition. We continue to verify that $f$ is a rerooting with respect to $S$.

Consider $a\in\st\setminus S$ and $b\in S$ with $a<b$. We claim that $f(a) \perp f(b)$. Pick $c\in S$ with $c>b$. Then $f({c})<f(b)$ and $B(a,b,c)$ imply that $f(a)>f(b)$ or $f(a)\perp f(b)$. The first case is impossible by the definition of $S$, and so $f(a)\perp f(b)$, verifying the claim. 
Next, consider $a\in\st\setminus S$ and $b\in S$ with $a\perp b$. Picking $c\in S$ with $B(a,c,b)$, we derive that $f(a)< f(b)$.

Let $x,y\in \st\setminus S$ with $x<y$. Pick $z\in S$ such that $y<z$. 
Note that $B(x,y,z)$ and following the claim above $f(x) \perp f(z)$ and $f(y)\perp f(z)$. 
Then $B(f(x),f(y),f(z))$, $f(x)\perp f(z)$, and $f(y)\perp f(z)$ imply that $f(x)<f(y)$. Given $x,y\in \st\setminus S$ with $x\perp y$, we can pick $z\in S$ such that $x<z$ and $y<z$. Then 
following the claim above and knowing that $f$ preserves $B$ and $\neg B$ (note Lemma~\ref{lem:not-b}), we have  
\[f (x) \perp f(z), f(y) \perp f (z), \neg B(f(x), f(y), f(z)), \text{ and } \neg B(f(y), f(x), f(z)) \] 
that together imply $f(x) \perp f(y)$.
\end{proof}

\begin{cor}\label{cor:rerooting}
$\Aut(\st;B)$ consists precisely of the surjective rerootings.
\end{cor}
\begin{proof} 
Every element of $\Aut(\st;B)$ is a rerooting 
by Proposition~\ref{prop:rerooting} and is surjective. Conversely, let $\alpha$ be a surjective rerooting with respect to 
$S$. Let $\beta$ be a rerooting with respect to
$\alpha[S]$. Then 
$\beta \circ \alpha[\st]$ is isomorphic to $\st$, so $\beta$ can be chosen surjectively. Since
$\beta \circ \alpha$ is an automorphism of 
$(\st;B)$, there is $\gamma \in \Aut(\st;B)$ such that $\gamma \circ \beta \circ 
\alpha$ is the identity, so $\alpha$ has
the inverse
$\gamma \circ \beta \in \Emb(\st;B)$ and thus is an automorphism of $(\st;B)$.
\end{proof}

\begin{cor}\label{cor:rerootinggenerate}
Let $e \in \Emb(\st;B)$ be such that
it does not preserve $<$. Then
$\{e\} \cup \Aut(\st;<)$ generates
$\Emb(\st;B)$.
\end{cor}
\begin{proof}
Let $e' \in \Emb(\st;B)$. 
By Proposition~\ref{prop:rerooting}, both $e$ and 
$e'$ are rerootings with respect to chains $S$ and $S'$. Then by the homogeneity of $(\st;\leq,C)$ 
for every finite subset $F$ of $\st$
there exists an automorphism $\alpha$
of $(\st;\leq,C)$ such that for 
every $x \in F$ it holds that
$x \in S'$ if and only if $\alpha(x) \in S$. 
By the definition of the rerooting
operation and again by homogeneity of 
$(\st;\leq,C)$
there exists an automorphism $\beta$ of $(\st;\leq,C)$ such that
$e'(x) = \beta(e(\alpha(x)))$ for all $x \in F$. Hence, by topological closure, $\{e\} \cup \Aut(\st;<)$ generates $e'$.
\end{proof}

\begin{corollary}\label{cor:endb}
$\End(\mS_2; B)=\Emb(\mS_2; B)=\overline{\Aut(\mS_2; B)}$.
\end{corollary}
\begin{proof}
Lemma~\ref{lem:not-b} shows that 
$\End(\mS_2;B) = \Emb(\mS_2;B)$. 
From Propositions~\ref{prop:rerooting} and~\ref{cor:rerooting} it follows that the restriction of any self-embedding of $(\mS_2; B)$ to a finite subset of $\st$ extends to an automorphism, and hence $\Emb(\mS_2; B)=\overline{\Aut(\mS_2; B)}$ by the definition of the pointwise convergence 
topology. 
\end{proof}

\subsection{Ramsey-theoretic analysis}
\label{sect:ramsey-analysis}

\subsubsection{Canonical functions without constants} Every canonical function $f\colon (\mS_2; \leq, C, \prec)\rightarrow (\mS_2; \leq, C,\prec)$ induces a function on the 3-types of $(\mS_2; \leq, C,\prec)$. Our first lemma shows that only few functions on those 3-types are induced by canonical functions, i.e., there are only few behaviors of canonical functions.

\begin{defn}
We call a function $f\colon \mS_2\rightarrow \mS_2$
\begin{itemize}
\item \emph{flat} \dff its image induces an antichain in $(\mS_2; \leq)$;
\item \emph{thin} \dff its image induces a chain in $(\mS_2; \leq)$.
\end{itemize}
\end{defn}

\begin{lemma}\label{lem:nocnst}
Let $f\colon (\mS_2; \leq, C, \prec)\rightarrow (\mS_2; \leq, C,\prec)$ be an injective canonical function. Then either $f$ is flat, or $f$ is thin, or $f\in\End(\mS_2;<,\perp)$.
\end{lemma}
\begin{proof}
Let $u_1,u_2,v_1,v_2\in\mS_2$ be so that $u_1<u_2$, $v_1\perp v_2$, and $v_1\prec v_2$. 
By the homogeneity of $(\st;\leq,C,\prec)$ all pairs
of distinct elements have the same
type as $(u_1,u_2),(v_1,v_2),(v_2,v_1)$, or $(u_2,u_1)$, and by canonicity pairs of equal type are sent
to pairs of equal type.
Hence, 
if $f(u_1)\perp f(u_2)$ and $f(v_1)\perp f(v_2)$, then $f$ is flat by canonicity. 
If $f(u_1)\nparallel f(u_2)$ and $f(v_1)\nparallel f(v_2)$, then $f$ is thin. It remains to check the following cases.

{\emph {Case 1: $f(u_1)\perp f(u_2)$ and $f(v_1)< f(v_2)$.}} Let $x,y,z\in \mS_2$ be such that $x<y$, $x\perp z$, $y\perp z$, $z\prec x$, and $z\prec y$. Then $f(x)\perp f(y)$, $f(x)>f(z)$, and $f(y)> f(z)$, in contradiction with the axioms of the semilinear order.

{\emph {Case 2: $f(u_1)\perp f(u_2)$ and $f(v_1)> f(v_2)$.}} Let $x,y,z\in \mS_2$ be such that $x<y$, $x\perp z$, $y\perp z$, $x\prec z$, and $y\prec z$. Then $f(x)\perp f(y)$, $f(x)>f(z)$, and $f(y)>f(z)$, in contradiction with the axioms of the semilinear order.

{\emph {Case 3: $f(u_1)<f(u_2)$ and $f(v_1)\perp f(v_2)$.}} Then $f$ preserves $<$ and $\perp$.

{\emph {Case 4: $f(u_1)>f(u_2)$ and $f(v_1)\perp f(v_2)$.}} Let $x,y,z\in \mS_2$ such that $x\perp y$, $x\prec y$, $x<z$, and $y< z$. Then $f(x)\perp f(y)$,  $f(x)>f(z)$, and $f(y)> f(z)$, in contradiction with the axioms of the semilinear order.
\end{proof}

\subsubsection{Canonical functions with constants}

\begin{lem}\label{lem:generateflat}
Let $f\colon \mS_2\rightarrow \mS_2$ be a function. If $f$ preserves incomparability but not comparability in $(\mS_2;\leq)$, then $\{f\}\cup\Aut(\st;\leq)$ generates a flat function. If $f$ preserves comparability but not incomparability in $(\mS_2;\leq)$, then $\{f\}\cup\Aut(\st;\leq)$ generates a thin function.
\end{lem}
\begin{proof}
We show the first statement; the proof of the second statement is analogous. 
We first claim that for any finite set $A\subseteq\st$, $\{f\} \cup \Aut(\st;\leq)$ generates a function which sends $A$ to an antichain. To see this, let $A$ be given, and pick $a,b\in\st$ such that $a<b$ and $f(a)\perp f(b)$. If $A$ contains elements $u,v$ with $u<v$, then there exists $\alpha\in\Aut(\st;\leq)$ so that $\alpha(u)=a$ and $\alpha(v)=b$ since the map that sends $(u,v)$ to $(a,b)$ is an isomorphism between substructures of the homogeneous structure $(\st;\leq,C)$. The function $f\circ \alpha$ sends $A$ to a set which has fewer pairs $(u,v)$ satisfying $u<v$ than $A$. Repeating this procedure on the image of $A$, and so forth, and composing functions we obtain a function which sends $A$ to an antichain. 

Now let $\{s_0,s_1,\ldots\}$ be an enumeration of $\st$, and pick for every $n\geq 0$ a function $g_n$ generated by $\{f\}\cup\Aut(\st;\leq)$ which sends $\{s_0,\ldots,s_n\}$ to an antichain. Since $(\mS;\leq)$ 
has finitely many orbits of $n$-tuples,
an easy consequence of K\"onig's
tree lemma shows that 
we may assume that for all $n\geq 0$ and all $i,j\geq n$ the type of the tuple $(g_i(s_0),\ldots,g_i(s_n))$ equals the type of $(g_j(s_0),\ldots,g_j(s_n))$ in $(\mS;\leq)$. By composing with automorphisms of $(\mS;\leq)$ from the left, we may even assume that these tuples are equal. But then the sequence $(g_n)_{n\in\omega}$ converges to a flat function.
\end{proof}

\begin{defn}
When $n\geq 1$ and $R\subseteq \mS_2^n$ is an $n$-ary relation, then we say that \emph{$R(X_1,\ldots,X_n)$ holds for sets $X_1,\ldots,X_n\subseteq \mS_2$} \dff $R(x_1,\ldots,x_n)$ holds whenever $x_i\in X_i$ for all $1\leq i\leq n$. We also use this notation when some of the $X_i$ are elements of $\mS_2$ rather than subsets, in which case we treat them as singleton subsets.
\end{defn}

\begin{defn}
For $a\in \mS_2$, we set
\begin{itemize}
\item $U_<^a:=\{p\in \mS_2\mid p<a\}$;
\item $U_>^a:=\{p\in \mS_2\mid p>a\}$;
\item $U_{\perp,\prec}^a:=\{p\in \mS_2\mid p\perp a \wedge p\prec a\}$;
\item $U_{\perp,\succ}^a:=\{p\in \mS_2\mid  p\perp a \wedge a\prec p\}$;
\item $U_{\perp}^a:=U_{\perp,\succ}^a\cup U_{\perp,\prec}^a$.
\end{itemize}
\end{defn}

The first four sets defined above are precisely the infinite orbits of $\Aut(\mS_2; \leq,\prec, a)$.

\begin{lemma}\label{lem:onecnst}
Let $a\in \mS_2$, and let $f\colon (\mS_2; \leq, C, \prec, a)\rightarrow (\mS_2; \leq, C,\prec)$ be an injective canonical function. Then one of the following holds:
\begin{enumerate}
\item $\{f\}\cup \Aut(\mS_2; \leq)$ generates a flat or a thin function;
\item $f\in\End(\mS_2;<,\perp)$;
\item $f(a)\not < f[U_{>}^a]$ and $f\rest_{\st\setminus\{a\}}$ behaves like a rerooting function with respect to $U_>^a$ in the
following sense: whenever $g$ is such a rerooting function, and $F$ is finite, then
there exists $\alpha \in \Aut(\mS_2; \leq)$ such that $\alpha g \rest_F=f \rest_F$.
\end{enumerate}
Moreover, if $f(a)\not > f[U_<^a]$ and $f(a)\not > f[U_{>}^a]$, then $\{f\}\cup \Aut(\mS_2; \leq)$ generates a flat or a thin function.
\end{lemma}
\begin{proof}
The structure induced by $U_<^a$ in $(\mS_2; \leq, C, \prec)$ is isomorphic to $(\mS_2; \leq, C, \prec)$ 
($U_<^a$ induces in $(\mS_2; \leq)$ a dense, unbounded, binary branching, and nice semilinear order without joins). 
The restriction of $f$ to this copy is canonical.
Since tuples on $U^a_<$ of equal type in
$(\mS_2;\leq,C,\prec)$ have equal type in 
$(\mS_2;\leq,C,\prec,a)$, such tuples
are sent to tuples of equal type 
in $(\mS_2;\leq,C,\prec)$ under $f$, by its canonicity. 
Picking a self-embedding $g$ of $(\mS_2;<,\prec)$ 
with image $U^a_<$, the composite function $fg$ is canonical from $(\mS_2;\leq,C,\prec)$
to $(\mS_2;\leq,C,\prec)$. 
If $f$ violates $<$ or $\perp$ on $U^a_<$, 
then $fg$ violates $<$ or $\perp$, 
and hence generates a flat or thin function 
by Lemma~\ref{lem:nocnst}.
Hence, we may assume that 
$f$ preserves $<$ and $\perp$ on $U_<^a$.  

When $u,v\in U_{\perp,\prec}^a$ satisfy $u<v$, then there exists a subset of $U_{\perp,\prec}^a$ containing $u$ and $v$ which induces an isomorphic copy of $(\mS_2; \leq, C, \prec)$. As above, we may assume that $f$ preserves $<$ and $\perp$ on this subset, and hence $f(u)<f(v)$. If $u,v\in U_{\perp,\prec}^a$ satisfy $u\perp v$, then there exist subsets $R,S$ of $U_{\perp,\prec}^a$ containing $u$ and $v$, respectively,  such that both $R$ and $S$ induce isomorphic copies of $(\mS_2; \leq, C, \prec)$ and such that for all $r\in R$ and $s\in S$ the type of $(r,s)$ equals the type of $(u,v)$ in $(\mS_2; \leq, C, \prec)$. Assuming as above that $f$ preserves $<$ and $\perp$ on both copies, $f(u)<f(v)$ would imply $f[R]<f[S]$, which is in contradiction with the axioms of a semilinear order. Hence, we may assume that $f$ preserves $<$ and $\perp$ on $U_{\perp,\prec}^a$, and by a similar argument also on $U_{\perp,\succ}^a$.

The sets $U_{\perp,\prec}^a$, $U_{\perp,\succ}^a$, and $U_<^a$ are pairwise incomparable, and $f$ cannot violate the relation $\perp$ between them, since by the canonicity of $f$ this would contradict the axioms of the semilinear order. Thus we may assume that $f$ preserves $<$ and $\perp$ on $U_{\perp}^a \cup U_<^a$. Moreover, for no $p\in\{a\}\cup U_{>}^a$ we have $f({p})<f[U_{\perp,\prec}^a]$, $f({p})<f[U_{\perp,\succ}^a]$, or $f({p})<f[U_<^a]$, again by the properties of semilinear orders.

Assume that $U_{>}^a$ is mapped to an antichain by $f$. Then the canonicity of $f$ implies that $f[U_{>}^a]\perp f[U_{\perp}^a\cup U_<^a]$, as all other possibilities are in contradiction with the axioms of the semilinear order. In particular, $f$ then preserves $\perp$ on $\st\setminus\{a\}$.
Given a finite $A\subseteq \mS_2$ which is not an antichain, there exists $\alpha\in\Aut(\mS_2; \leq)$ such that $\alpha[A]\subseteq \mS_2\setminus\{a\}$ and two comparable points of $A$ are mapped into $U_{>}^a$ by $\alpha$. Thus $f\circ \alpha$ preserves $\perp$ on $A$, and it maps at least one comparable pair in $A$ to an incomparable one. As in Lemma~\ref{lem:generateflat}, we see that $\{f\}\cup \Aut(\mS_2; \leq)$ generates a flat function. So we may assume that the order on $U_{>}^a$ is either preserved or reversed by $f$. The rest of the proof is an analysis of the possible behaviours of $f$ in these two cases. In order to talk about the behaviour of $f$, we choose elements $u_1\in U_{\perp,\prec}^a$, $u_2\in U_{\perp,\succ}^a$ and  $z_1, z_2\in U_{>}^a$ such that $z_1<z_2$, $u_i\perp z_1$, and $u_i< z_2$ for $i\in\{1,2\}$; see Figure~\ref{fig:onecnst}. 

\begin{figure}
\begin{center}
\includegraphics[scale=.5]{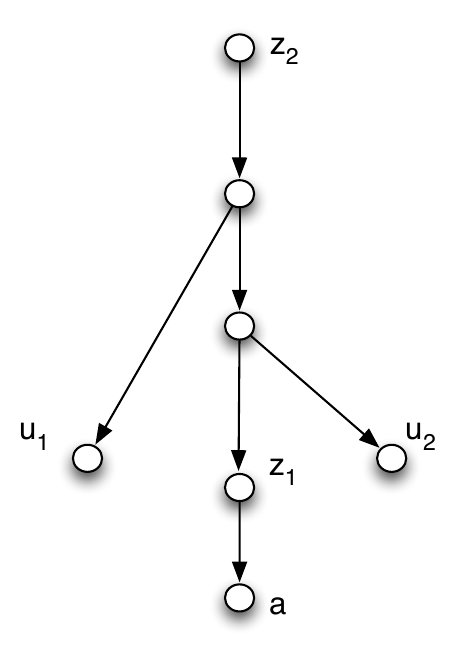}
\end{center}
\caption{Illustration for the case distinction in the proof of Lemma~\ref{lem:onecnst}.}
\label{fig:onecnst}
\end{figure}

{\emph {Case 1: $f$ preserves the order on $U_{>}^a$.}} If $f(u_1)<f(z_1)$, then by transitivity of $<$ and canonicity of $f$ we have that $f[U_{\perp,\prec}^a]<f[U_{>}^a]$.   
Given a finite $A\subseteq \mS_2$ which is not a chain, there exists $\alpha\in\Aut(\mS_2; \leq)$ such that $\alpha[A]\subseteq U_{\perp,\prec}^a\cup U_{>}^a$ and such that $\alpha(x)\in U_{\perp,\prec}^a$ and $\alpha(y)\in U_{>}^a$ for some elements $x,y\in A$ with $x\perp y$. 
Thus $f\circ \alpha$ preserves $<$ on $A$, and it maps at least one incomparable pair in $A$ to a comparable one. As in Lemma~\ref{lem:generateflat}, we conclude that $\{f\}\cup \Aut(\mS_2; \leq)$ generates a thin function. We can argue similarly when $f(u_2)<f(z_1)$. Thus we may assume that $f(u_i)\perp f(z_1)$ for $i\in \{1,2\}$. If $f(u_i)\perp f(z_2)$ for some $i\in \{1,2\}$, then a similar argument shows that $\{f\}\cup \Aut(\mS_2; \leq)$ generates a flat function. 
Hence, we may assume that $f(u_i)<f(z_2)$ for $i\in \{1,2\}$, and so $f$ preserves $<$ and $\perp$ on $U_\perp^a\cup U_>^a$.

Assume that $f[U_<^a]\perp f[U_{>}^a]$. Given a finite $A\subseteq \mS_2$ which is not an antichain, there exists $\alpha\in\Aut(\mS_2; \leq)$ such that $\alpha[A]\subseteq \mS_2\setminus\{a\}$ and such that $\alpha(x)\in U_<^a$ and $\alpha(y)\in U_{>}^a$ for some $x,y\in A$ with $x<y$. 
Thus $f\circ \alpha$ preserves $\perp$ on $A$, and it maps at least one comparable pair in $A$ to an incomparable one. The proof of Lemma~\ref{lem:generateflat} shows that $\{f\}\cup \Aut(\mS_2; \leq)$ generates a flat function. 
So we may assume that $f[U_<^a]< f[U_{>}^a]$, and consequently, $f$ preserves $<$ and $\perp$ on $\mS_2\setminus \{a\}$.

If $f(a)>f[U_{>}^a]$, then by transitivity of $<$ we have $f(a)>f[\mS_2\setminus \{a\}]$, and we can easily show that $\{f\}\cup \Aut(\mS_2; \leq)$ generates a thin function. 
Similarly, if $f(a)\perp f[U_{>}^a]$, then by the axioms of the semilinear order we have $f(a)\perp f[\mS_2\setminus \{a\}]$, and $\{f\}\cup \Aut(\mS_2; \leq)$ generates a flat function. 
Thus we may assume that $f(a)<f[U_{>}^a]$. 
If $f(a)> f[U_{\perp,\prec}^a]$ or $f(a)> f[U_{\perp,\succ}^a]$, then by transitivity of $<$ we have $f[U_{\perp,\prec}^a]<f[U_{>}^a]$ or $f[U_{\perp,\succ}^a]<f[U_{>}^a]$, a contradiction. Hence, $f(a)\perp f[U_{\perp}^a]$. 
Finally, if $f(a)\perp f[U_<^a]$, then $\{f\}\cup \Aut(\mS_2; \leq)$ generates a flat function. 
Thus we may assume that $f(a)>f[U_<^a]$, and so $f$ preserves $<$ and $\perp$, proving the lemma. 

{\emph {Case 2: $f$ reverses the order on $U_{>}^a$.}} If $f(u_1)\perp f(z_1)$, then by $f(z_2)<f(z_1)$ and the axioms of the semilinear order we have that $f(u_1)\perp f(z_2)$. 
Moreover, $f\rest_{U_{\perp,\prec}^a \cup U_{>}^a}$ preserves $\perp$. Since the comparable elements $u_1, z_2$ are sent to incomparable ones, the standard iterative argument shows that $\{f\}\cup \Aut(\mS_2; \leq)$ generates a flat function. An analogous argument works if $f(u_2)\perp f(z_1)$. 
Thus we may assume that $f(u_i)<f(z_1)$ for $i\in \{1,2\}$. If $f(u_i)<f(z_2)$ for some $i\in \{1,2\}$, then a similar argument shows that $\{f\}\cup \Aut(\mS_2; \leq)$ generates a thin function. 
Thus we may assume that $f(u_i)\perp f(z_2)$ for $i\in \{1,2\}$, and $f\rest_{U_\perp^a\cup U_>^a}$ behaves like a rerooting.

Assume that $f[U_<^a]<f[U_{>}^a]$. Let $A\subseteq \mS_2$ be finite. Pick a minimal element $b \in A$, and let $C\subseteq A$ be those elements $c\in A$ with $b\leq c$. 
Let $\alpha\in\Aut(\mS_2; \leq)$ be such that $\alpha(b)\in U_<^a$, $\alpha[C\setminus\{b\}]\subseteq U_{>}^a$ and $\alpha[A\setminus C]\subseteq U_{\perp}^a$. 
Then there exists $\beta\in\Aut(\mS_2; \leq)$ such that $\beta\circ f\circ \alpha[C]\subseteq U_{>}^a$ and $\beta\circ f\circ \alpha[A\setminus C]\subseteq U_{\perp}^a$. 
Let $g:=f\circ \beta\circ f\circ \alpha$. 
Then $g\rest_{A\setminus\{b\}}$ preserves $<$ and $\perp$, and $g(b)\geq g[A]$. 
By iterating such steps, $A$ can be mapped to a chain. Hence, as in Lemma~\ref{lem:generateflat}, $\{f\}\cup \Aut(\mS_2; \leq)$ generates a thin function. Thus we may assume that $f[U_<^a]\perp f[U_{>}^a]$. 
By replacing $U_<^a$ with $\{a\}$ in this argument, one can show that if $f(a)<f[U_{>}^a]$, then $\{f\}\cup \Aut(\mS_2; \leq)$ generates a thin function. Thus we may assume that $f(a)\not < f[U_{>}^a]$, and so Item~(3) applies.

To show the second part of the lemma, suppose that $f(a)\not > f[U_<^a]$ and $f(a)\not > f[U_{>}^a]$. 
Then $f$ violates $<$, thus Item (2) cannot hold for $f$.  Hence, either $\{f\}\cup \Aut(\mS_2; \leq)$ generates a flat or a thin function, or the conditions in Item~(3) hold for $f$. We assume the latter. In particular, $f(a)\perp f[U_\perp^a]$, by the axioms of the semilinear order, and hence 
$f(a)\perp f[U_{>}^a]$.

Let $A\subseteq \mS_2$ be finite such that $A$ is not an antichain. Pick some $x\in A$ with is maximal in $A$ with respect to $\leq$ and such that there exists $y\in A$ with $y<x$. 
Let $\alpha\in\Aut(\mS_2; \leq)$ be such that $\alpha(x)=a$. Then $f\circ \alpha$ preserves $\perp$ on $A$, and $f(y)\perp f(x)$. Hence, iterating such steps $A$ can be mapped to an antichain, and $\{f\}\cup \Aut(\mS_2; \leq)$ generates a flat function.
\end{proof}

\subsubsection{Applying canonicity}

\begin{lemma}\label{lem:climb1end}
Let $f\colon\mS_2\rightarrow \mS_2$ be an injective function that violates $<$. Then either $\{f\}\cup \Aut(\mS_2; \leq)$ generates a flat or a thin function, or $\{f\}\cup \Aut(\mS_2; \leq)$ generates $\End(\mS_2; B)$.
\end{lemma}
\begin{proof}
It is easy to see that if $f$ preserves comparability and incomparability, then $f$ cannot violate $<$. If $f$ preserves comparability and violates incomparability, then $\{f\}\cup \Aut(\mS_2; \leq)$ generates a thin function by Lemma~\ref{lem:generateflat}. 
Thus we may assume that $f$ violates comparability. 
Let $a,b\in\st$ such that $a<b$ and $f(a)\perp f(b)$. 
According to Proposition~\ref{prop:canfcts}, there exists a canonical function $g\colon(\mS_2; \leq, C, \prec, a, b)\rightarrow (\mS_2; \leq, C,\prec)$ that is generated by $\{f\}\cup \Aut(\mS_2; \leq)$ such that $g(a)\perp g(b)$. 
The set $U_<^b$ induces in $(\mS_2; \leq, C, \prec, a)$ a structure that is isomorphic to $(\mS_2; \leq, C, \prec, a)$, and the restriction of $g$ to this set is canonical. 
By Lemma~\ref{lem:onecnst} either $\{g\}\cup \Aut(\mS_2; \leq)$ generates a thin or a flat function (case (1)), or a rerooting (case (3)), in which case 
$\{g\} \cup \Aut(\mS_2;\leq)$ generates $\End(\mS_2;B)$ 
by Proposition~\ref{prop:rerooting}, Corollary~\ref{cor:rerootinggenerate}, and Corollary~\ref{cor:endb}, or $g$ preserves $<$ and $\perp$ on $U_<^b$ (case (2)). In the first two cases we are done
so we may assume the latter. By a similar argument, either $\{g\}\cup \Aut(\mS_2; \leq)$ generates a thin or a flat function, or a rerooting, or $g$ preserves $<$ and $\perp$ on $U_<^a\cup U_\perp^b\cup U_>^b\cup\{b\}$.
However, the latter is impossible as it would imply that $g(t)<g(a)$ and $g(t)< g(b)$ for all $t\in U_<^a$ while $g(a)\perp g(b)$, which is in contradiction with the axioms of the semilinear order.
\end{proof}

Next, we study injective functions $f$ that violate $B$. 
The main result will be Lemma~\ref{lem:climb2end} 
stating that such functions generate flat or thin functions. 
The following fact is a special case of \cite{AdelekeNeumann}, Corollary 20.7;
we just sketch an argument in the present terminology for the benefit of the reader. 
\begin{prop}
\label{prop:B3hom}
 Every isomorphism between 3-element substructures of $(\st;B)$ extends to an automorphism of $(\st;B)$.
\end{prop}
\begin{proof}
Take 3-element sets $D,D' \subseteq \st$ which induce isomorphic structures in $(\st;B)$ and $p \colon D \rightarrow D'$ an isomorphism between them. Whatever the isomorphism types induced by $D,D'$ are in $(\st;\leq,C)$ it can be checked that one can apply a surjective rerooting $h$ to $D$ so that the structures induced by $h[D]$ and $D'$ in $(\st;\leq,C)$ are isomorphic. By Corollary \ref{cor:rerooting} this $h$ is in $\Aut(\st;B)$. It follows from the homogeneity of $(\st;\leq,C)$ that there 
exists a $\beta \in \Aut(\st;\leq,C)$ such that $\beta \circ h [D] = D'$, hence $\beta \circ h$ is an automorphism of $(\st;B)$ extending $p$. 
\end{proof}
To ease notation, define ternary relations $K$ and $L$ on $\st$ by 
\begin{eqnarray*}
K(x,y,z) &\Leftrightarrow &x \neq y \neq z \neq x \wedge \neg B(x,y,z) \wedge \neg B(y,z,x) \wedge \neg B(z,x,y);\\
L(x,y,z) &\Leftrightarrow &B(x,y,z) \vee B(y,z,x) \vee B(z,x,y). 
\end{eqnarray*}
Note that $K$ and $L$ are mutually exclusive and for any distinct $x,y,z$ in $\st$, we have $K(x,y,z) \vee L(x,y,z)$. 
\begin{lemma}
\label{lem:killKthin}
Suppose that $f \colon \st \rightarrow \st$ is an injective function such that there are $a,b,c \in \st$ distinct with $K(a,b,c)$ and $B(f(a),f(b),f(c))$, yet there is no $r,s,t \in \st$ with $B(r,s,t)$ and $K(f(r),f(s),f(t))$. Then $\{f\} \cup \Aut(\st;B)$ generates a thin function.
\end{lemma}
\begin{proof}
It suffices to prove that for any finite set $A \subseteq \st$ of cardinality at least three, $\{f\} \cup \Aut(\st;B)$ generates a function $g$ such that $g[A]$ is a chain in $(\st; <)$; by the $\omega$-categoricity of $(\st;B)$, one can then apply K\"{o}nig's tree lemma, as in the second paragraph of the proof of Lemma \ref{lem:generateflat}, to prove that $\{f\} \cup \Aut(\st;B)$ generates a thin function.

Fix a finite set $A \subseteq \st$ of cardinality at least 3.  
Suppose that $u,v,w \in A$ are such that $K(u,v,w)$. By Proposition \ref{prop:B3hom} there is a $\beta \in \Aut(\st;B)$ such that $(\beta(u),\beta(v),\beta(w)) = (a,b,c)$ and so $B(f \circ \beta(u), f \circ \beta(v), f \circ \beta(w))$. Hence $f \circ \beta [A]$ contains fewer triples $x,y,z$ satisfying $K(x,y,z)$ than $A$ does. Iterating this procedure on $f \circ \beta [A]$ and so forth a finite number of times, composing functions as we go, we obtain a function generated by $\{f\} \cup \Aut(\st;B)$ which sends $A$ to a set $A'$ such that for all distinct $u,v,w$ in $A'$ we have $L(u,v,w)$. If there is no $u,v,w \in A$ such that $K(u,v,w)$, then we have $L(u,v,w)$ for all $u,v,w$ in $A$, so we may set $A'$ to be the image of $A$ under the identity. In either case, clearly $A' = C_1 \cup C_2$ is the disjoint union of at most two mutually incomparable chains $C_1$ and $C_2$. Then there is a surjective rerooting $\alpha$ with respect to a maximal chain in $(\st;<)$ extending $C_2$; by Corollary \ref{cor:rerooting} this $\alpha$ is in $\Aut(\st;B)$. We then have that $\alpha[C_1] < \alpha[C_2]$, so $\alpha[A']$ is a chain and we are done. 
\end{proof}

\begin{lemma}\label{lem:climb2end}
Let $f\colon\mS_2\rightarrow \mS_2$ be an injective function that violates $B$. Then $\{f\}\cup\Aut(\mS_2; B)$ generates a flat or a thin function.
\end{lemma}
\begin{proof}
Suppose there are no $r,s,t \in \st$ such that $B(r,s,t)$ and $K(f(r),f(s),f(t))$. As $f$ violates $B$, we may assume there are some $a,b,c \in \st$ such that $B(a,b,c)$ and $B(f(b),f(a),f(c))$. By Proposition \ref{prop:B3hom}, there are $\alpha, \beta \in \Aut(\st;B)$ so that $\alpha(a)<\alpha(b)<\alpha(c)$ and $\beta(f(b))<\beta(f(a))<\beta(f(c))$. Hence we may assume that $a<b<c$ and $f(b)<f(a)<f(c)$, if necessary replacing $a,b,c$ by $\alpha^{-1}(a),\alpha^{-1}(b),\alpha^{-1}(c)$ and $f$ by $\beta \circ f \circ \alpha^{-1}$ and relabeling. Take $d$ such that $a \perp d < b$ so that $K(a,b,d)$ and consider where $f$ takes $d$. If $f(d) < f(b)$ then we have $B(f(d),f(b),f(a))$, so applying Lemma \ref{lem:killKthin} we know that a thin function is generated. If $f(d) > f(a)$ or $f(d) \perp f(a)$, then $B(f(d),f(a),f(b))$ and a thin function is generated by Lemma \ref{lem:killKthin}. If $f(b) < f(d) < f(a)$, then $B(f(b),f(d),f(a))$ and a thin function is generated. By semilinearity the remaining possibility is $f(b) \perp f(d) < f(a)$. Again by semilinearity, we have that $a \perp d < c$ and $f(d) < f(a) < f(c)$. In particular we have that $K(d,a,c)$ and $B(f(d),f(a),f(c)$ so, by Lemma \ref{lem:killKthin}, $\{f\} \cup \Aut(\st;B)$ generates a thin function.

We may now assume there are 
$a,b,c\in \mS_2$ be such that $B(a,b,c)$ and $K(f(a),f(b),f(c))$. 
It follows from Proposition \ref{prop:B3hom} that there exist $\alpha, \beta \in \Aut(\st;B)$ such that $\alpha(a) < \alpha(b) < \alpha(c)$ and that $\{\beta(f(a)),\beta(f(b)),\beta(f(c))\}$ induces an antichain. Replacing $f$
 by $\beta\circ f\circ\alpha^{-1}$, we may assume that there are $a,b,c\in \mS_2$ such that $a<b<c$ and such that $\{f(a), f(b), f(c)\}$ induces an antichain. By Proposition~\ref{prop:canfcts}, there is a canonical function $g\colon (\mS_2; \leq, C, \prec, a, b, c)\rightarrow (\mS_2; \leq, C)$ that is generated by $\{f\}\cup \Aut(\mS_2; \leq)$ such that $\{g(a), g(b),g(c)\}$ induces an antichain.
 
By the axioms of the semilinear order, at most one $y\in\{g(a), g(b), g({c})\}$ can satisfy $y>g[U_<^a]$ and at most one such element can satisfy $y>g[U_>^c]$. 
Hence, there exists an $x\in \{a, b, c\}$ such that $g(x)\not > g[U_<^a]$ and $g(x)\not> g[U_>^c]$. 
The set $X:=U_<^a\cup \{x\}\cup U_>^c\cup U_{\perp}^c$ induces in $(\mS_2; \leq)$ a structure isomorphic to $(\mS_2; \leq)$ because is a dense, unbounded, binary branching, nice, semilinear order without joins, and consequently 
$X$ induces in $(\mS_2; \leq,C,\prec)$
a structure isomorphic to $(\mS_2;\leq,C,\prec)$. Moreover, $g\rest_{X}$ is canonical as a function from $(\mS_2; \leq, C, \prec,x)$ to $(\mS_2; \leq, C, \prec)$; this can be
shown by the same line of argument as given
in the first paragraph of Lemma~\ref{lem:onecnst}.
According to the second part of Lemma~\ref{lem:onecnst}, $\{g\}\cup \Aut(\mS_2; \leq)$ generates a flat or a thin function.
\end{proof}


\subsection{Endomorphisms and the proof of Theorem~\ref{thm:4case}}
\label{sect:4case}

\begin{prop}\label{prop:embeddings}
Let $\Gamma$ be a reduct of $(\mS_2; \leq)$. Then one of the following holds.
\begin{enumerate}
\item $\End(\Gamma)$ contains a flat or a thin function.
\item $\End(\Gamma)= \overline{\Aut(\mS_2; \leq)}$.
\item $\End(\Gamma)= \overline{\Aut(\mS_2; B)}$.
\end{enumerate}
\end{prop}
\begin{proof}
Assume that there exist $x,y\in\st$ with $x<y$ and $f\in \End(\Gamma)$ such that $f(x)=f(y)$. By collapsing comparable pairs one-by-one using $f$ and automorphisms of $(\mS_2; \leq)$, 
it is possible to obtain for every finite
$F \subseteq \st$ a function $f' \in \End(\Gamma)$ such that $f'[F]$ induces an antichain. Hence, 
$\End(\Gamma)$ contains 
a flat function because $\End(\Gamma)$ is topologically closed in $(\st)^{\st}$. 

Similarly, if there exists a pair of elements $x\perp y$ and $f\in \End(\Gamma)$ such that $f(x)=f(y)$, then $\{f\}\cup \Aut(\mS_2; \leq)$ generates a thin function. Hence, we may assume that every endomorphism of $\Gamma$ is injective.
If $\End(\Gamma)$ preserves $<$ and $\perp$, then $$\End(\Gamma)= \Emb(\mS_2; \leq) = \Emb(\mS_2;\leq,C) = \overline{\Aut(\mS_2; \leq,C)} = \overline{\Aut(\mS_2; \leq)} \, .$$ If $\End(\Gamma)$ preserves $<$ and violates $\perp$, then $\End(\Gamma)$ contains a thin function. Thus we may assume that some $f\in \End(\Gamma)$ violates $<$. By Lemma~\ref{lem:climb1end} either $\End(\Gamma)$ contains a flat or a thin function, or $\Emb(\mS_2; B)\subseteq \End(\Gamma)$. Since $\Emb(\mS_2; B)= \overline{\Aut(\mS_2; B)}$, we may assume that $\Emb(\mS_2; B)\subsetneq \End(\Gamma)$, as otherwise Item (1) or (3) holds. Hence, there exists a function $f\in \End(\Gamma)$ that violates either $B$ or $\lnot B$. By Lemma~\ref{lem:not-b}, 
$f$ violates $B$, and then $\End(\Gamma)$ contains a flat or a thin function by Lemma~\ref{lem:climb2end}.
\end{proof}


\begin{lemma}\label{lem:canac}
Let $\Gamma$ be a reduct of $(\mS_2; \leq)$ which has a flat endomorphism. Then $\Gamma$ is homomorphically equivalent to a reduct of $(\mL_2; C)$.
\end{lemma}
\begin{proof}
Let $f$ be that endomorphism. By Zorn's lemma, there exists a maximal antichain $M$ in $\mS_2$ that contains the image of $f$. 
By definition $M$ induces in $(\mS_2;C)$ a structure $\Sigma$ which is isomorphic to $(\mL_2; C)$. 
The structure $\Delta$ with domain $M$ and all relations that are restrictions of the relations of $\Gamma$ to $M$ is a reduct of $\Sigma$, as $(\mS_2; \leq, C)$ has quantifier elimination (an $\omega$-categorical structure has quantifier-elimination if and only if it is homogeneous~\cite{Oligo}). 
The inclusion map of $M$ into $\mS_2$ is a homomorphism from $\Delta$ to $\Gamma$, and the function $f$ is a homomorphism from $\Gamma$ to $\Delta$.
\end{proof}

\begin{lemma}\label{lem:canch}
Let $\Gamma$ be a reduct of $(\mS_2; \leq)$ which has a thin endomorphism. Then $\Gamma$ is homomorphically equivalent to a reduct of the dense linear order.
\end{lemma}
\begin{proof}
Analogous to the proof of Lemma~\ref{lem:canac}, using the obvious fact that maximal chains in $(\mS_2; \leq)$ are isomorphic to $(\mathbb{Q}; \leq)$.
\end{proof}


\begin{proof}[Proof of Theorem~\ref{thm:4case}]
If $\End(\Gamma)$ contains a flat 
function, then $\Gamma$ is homomorphically equivalent to a reduct of $(\mL_2,C)$ by Lemma~\ref{lem:canac}; so we are in case (1). 
If $\End(\Gamma)$ contains a thin
function, then $\Gamma$ is homomorphically equivalent to a reduct of $(\mQ;\leq)$ by Lemma~\ref{lem:canch}; so we are in case (2). 
Otherwise, by Proposition~\ref{prop:embeddings},
we have that $\End(\Gamma) \in \{ \overline{\Aut(\mS_2; B)}, \overline{\Aut(\mS_2; \leq)}\}$,
and we are in case (3) and (4), respectively. In case (3), 
$\overline{\Aut(\mS_2;B)} = \End(\mS_2;B)$ by Corollary~\ref{cor:endb}, and hence $\Gamma$ is
existentially positively interdefinable
with $(\mS_2;B)$. 
In case (4), $$\overline{\Aut(\mS_2;\leq)} = \overline{\Aut(\mS_2;\leq,C)} =
\Emb(\mS_2;\leq,C) = \Emb(\mS_2;\leq) = \End(\st;<,\perp)$$
where the third equality holds because $C$ has an existential definition over $(\mS_2;\leq)$. 
Hence, $\Gamma$ is existentially positively interdefinable with $(\st;<,\perp)$. 
\end{proof}

\begin{proof}[Proof of Corollary~\ref{cor:mc-cores}]
We use Theorem~\ref{thm:4case}.
If $\Gamma$ is homomorphically equivalent to a reduct $\Gamma'$ of $({\mathbb Q}; <)$, then the model-complete core
of $\Gamma$ and $\Gamma'$ has only one element, or is interdefinable with one of the five structures from Cameron's theorem described earlier~\cite{tcsps-journal}. If $\Gamma$ is homomorphically equivalent to a reduct of $({\mathbb L}_2;C)$ then the model-complete core of $\Gamma$ and $\Gamma'$ has only one element, or is interdefinable with one of the three
structures $({\mathbb L};C)$, $({\mathbb L};D)$, or $({\mathbb L};=)$
by the mentioned result from~\cite{BodJonsPham}. Otherwise,
by Theorem~\ref{thm:4case},
$\End(\Gamma) \in 
\{\overline{\Aut(\st;B)},\overline{\Aut(\st; \leq)}\}$ and $\Gamma$ is its own 
model-complete core, and 
$\Gamma$ is existentially positively interdefinable with $(\st;B)$ or with 
$(\st;<,\perp)$. 
\end{proof}

\subsection{Embeddings and the proof of Theorem~\ref{thm:groups}}
\label{sect:groups}

\begin{lemma}\label{lem:thinemb}
Let $\Gamma$ be a reduct of $(\mS_2; \leq)$ with a thin self-embedding. Then $\Gamma$ is isomorphic to a reduct of $(\mathbb Q;<)$.
\end{lemma}
\begin{proof}
Let $f$ be a thin self-embedding of $\Gamma$. 
By Proposition~\ref{prop:canfcts},
$\{f\} \cup \Aut(\mS_2;\leq,\prec)$ generates a thin canonical function $g \colon (\mS_2; \leq, C, \prec)\rightarrow (\mS_2; \leq, C,\prec)$
which is also a self-embedding of $\Gamma$. 
There are four possible behaviours of $g$, as it can preserve or reverse $<$, and independently, it can preserve or reverse $\prec$ on incomparable pairs. 
In all four of these cases, the structure $\Sigma$ induced by the image of $g$ in $(\mS_2; \leq)$ is isomorphic to $(\mathbb Q;\leq)$. 
The structure $\Delta$ on this image whose relations are the restrictions of the relations of $\Gamma$ to $g[\st]$ is a reduct of $\Sigma$, as $(\mS_2; \leq, C)$ has quantifier elimination. The claim follows as $g$ is an isomorphism between $\Gamma$ and $\Delta$. 
\end{proof}

\begin{lemma}\label{lem:redq}
Let $\Gamma$ be a reduct of $(\mS_2; \leq)$ which is isomorphic to a reduct of $(\mathbb Q;<)$. Then $\Gamma$ is existentially interdefinable with $(\mS_2; =)$.
\end{lemma}
\begin{proof}
Pick any pairwise incomparable elements $a_1,\ldots,a_5\in \mS_2$. 
There exist $i,j \in \{1,\dots,5\}$ 
such that $C(a_ia_j,a_k)$ for all $k \in \{1,\dots,5\} \setminus \{i,j\}$. Then the mapping
 which flips $a_i,a_j$ and fixes the other three elements
 preserves $C$ and hence extends to an automorphism of $(\mS_2; \leq)$ by the
homogeneity of of $(\mS_2; \leq,C)$.
From Cameron's classification of the reducts of $(\mathbb Q;<)$ (\cite{Cameron5}, cf.~the description in Section~\ref{sect:results}) we know that the only automorphism group of such a reduct which can perform this is the full symmetric group, since all other groups fix at most one or all of five elements when they act on them. Hence, $\Aut(\Gamma)$ contains all permutations of $\mS_2$. Thus, all injections of $\mS_2$ are self-embeddings of $\Gamma$. In other words, $\Emb(\Gamma) =\Emb(\mS_2;=)$ and $\Gamma$ is existentially interdefinable with $(\mS_2; =)$.
\end{proof}

\begin{definition}
Let $R(x,y,z)$ be the ternary relation on $\mS_2$ defined by the formula 
$$
C(z,xy)\vee (x<z\wedge y<z)\vee (x\perp z\wedge y\perp z\wedge (x<y\vee y<x)).
$$
\end{definition}

\begin{proposition}\label{prop:notmc}
$(\st;R)$ and $(\st;\leq)$ are interdefinable. However, $(\st;R)$ is not model complete, i.e., it has a self-embedding which is not an element of $\overline{\Aut(\st;R)}$.
\end{proposition}
\begin{proof}
By definition, $R$ has a first-order definition in $(\st;\leq)$. To see the converse, observe that
 for $a,b\in\st$ we have that $a\leq b$ if and only if there exists no $c\in\st$ such that $R(b,c,a)$. Hence, $(\st;R)$ and $(\st;\leq)$ are interdefinable, and in particular, $\Aut(\st;R)=\Aut(\st;\leq)$.
 
To show that $(\st;R)$ is not model complete, let $f \colon \st \to \st$ be
such that $f[\st]$ is an antichain in $(\st;\leq)$ and such that $R(a,b,c)$ holds if and only if  $C(f({c}),f(a)f(b))$ for all $a,b,c\in\st$. 
The existence of such a function $f$ can be shown inductively as follows. We suppose that $f$
is already defined on a finite set $F$, and let
$x \in \mS_2 \setminus F$. The base case 
$F = \emptyset$ is trivial since $f(x)$ can be chosen arbitrarily in $\mS_2$. 
In the induction step, we distinguish the following cases. 
By inductive hypothesis and composing $f$ with itself if necessary,
we can assume that all the elements in $F$ are incomparable. 
\begin{itemize}
\item Case 1: $x$ is incomparable to all elements of $F$. In this case we are done.
\item Case 2: $x$ is comparable to some $y \in F$. Then $R(x,y,z)$ for all other $z \in F$ because all the elements in $F$ are incomparable. 
Define $f(x)$ such that $C(f(z),f(x)f(y))$ for all $z \in F$.
\end{itemize}

Clearly, $f$ is not an element of $\overline{\Aut(\st;R)}$, since it does not preserve comparability.
\end{proof}

%
The previous proposition is the reason for the special case concerning $R$ in the following lemma.

\begin{lemma}\label{lem:flatemb}
Let $\Gamma$ be a reduct of $(\mS_2; \leq)$ with a flat self-embedding. Then $\Gamma$ is isomorphic to a reduct of $(\mathbb Q;<)$, or it has a flat self-embedding that preserves $R$.
\end{lemma}
\begin{proof}
Let $f$ be the flat self-embedding. Proposition~\ref{prop:canfcts} implies that 
$\{f\} \cup \Aut(\mS_2; \leq,\prec)$ generates a function $g$ that
is canonical as a function from $(\mS_2; \leq, C,\prec)$ to $(\mS_2; \leq, C,\prec)$; since $g$ will also be a flat self-embedding of $\Gamma$,
we can replace $f$ by $g$ and can therefore assume
that already $f$ is canonical as a function from $(\mS_2; \leq, C,\prec)$ to $(\mS_2; \leq, C,\prec)$. 
Thus, $f$ either preserves $\prec$ between
any two incomparable elements, or it reverses $\prec$ between any two incomparable elements. 
A similar statement holds for pairs of comparable elements. Let $\alpha \in \Aut(\mS_2; \leq, C)$ be such that it reverses the order $\prec$ on incomparable pairs and preserves the order on comparable elements\footnote{Such an automorphism exists since the homogeneous structure $(\mS_2;\leq,C,\prec')$ where $x \prec' y \Leftrightarrow x \leq y \vee (x \perp y \wedge y \prec x)$ has the same age
as $(\mS_2;\leq,C,\prec)$; therefore 
these structures are isomorphic,
and any isomorphism gives an automorphism of $(\mS_2; \leq, C)$ as desired.}. Then $f \circ \alpha$ either preserves or reverses the order $\prec$. 
In the latter case, $\alpha \circ f$ preserves $\prec$ (as $f$ is flat), so in any case we may assume that $f$ preserves $\prec$.
To simplify notation, we shall write $x'$ instead of $f(x)$ for all $x\in \mS_2$, and we write $xy|z$ or $z|xy$ instead of $C(z,xy)$ for all $x,y,z\in\st$.

Let $a_1,\ldots,a_5\in \mS_2$ be so that $a_1\prec\cdots\prec a_5$ and so that $a_1\perp a_2$, $a_1,a_2<a_3$, $a_3\perp a_4$, and $a_1,\ldots,a_4<a_5$. See Figure~\ref{fig:lemma4-21-a1}, left side. 
\begin{figure}
\begin{center}
\includegraphics[scale=.5]{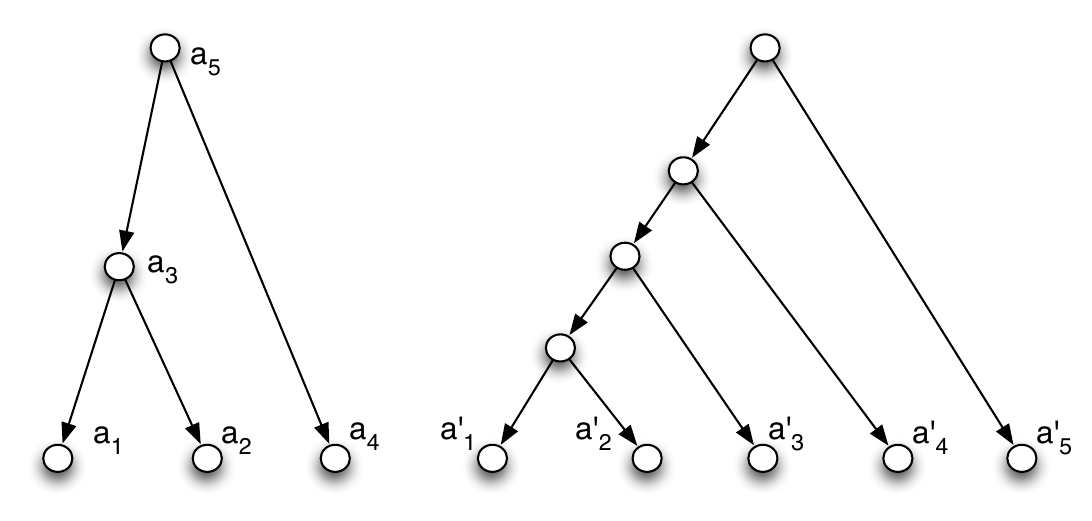}
\end{center}
\caption{First illustration for the analysis of the behaviour of the function $f$ from Lemma~\ref{lem:flatemb}; the picture on the right corresponds to case (A1).}
\label{fig:lemma4-21-a1}
\end{figure}
We shall analyse the possible behaviours of $f$ on these elements. 
Since $f$ preserves $\prec$ and by convexity of $\prec$, we have that either $a_1'a_2'|a_3'$ or $a_1'|a_2'a_3'$. 

We claim that in the first case, $a_2'a_3'|a_4'$. 
Pick  $x>a_2$ such that $a_1x|a_4$. 
Since $a_1'a_2'|a_3'$, we must have $a_1'a_2'|a_4'$ by the properties of $\prec$, and so $a_1'x'|a_4'$ by canonicity. 
Since $\prec$ extends $<$ and $f$ preserves $\prec$,
we have $a_1' \prec a_2' \prec x'$. These facts imply 
that $a_2'x'|a_4'$ by the properties of $\prec$, 
 and hence indeed $a_2'a_3'|a_4'$ by canonicity. This together with $a_1'a_2'|a_3'$ implies $a_1'a_3'|a_4'$. Since  $a_1'a_2'|a_3'$, we have $a_1'a_4'|a_5'$ by canonicity, leaving us with the following possibility which uniquely determines the type of the tuple $(a_1',\ldots,a_5')$ in $(\mS_2; \leq, C,\prec)$:
\begin{itemize}
\item[(A1)] $a_1'a_2'|a_3'$,  $a_1'a_3'|a_4'$, $a_1'a_4'|a_5'$;
\end{itemize}
see Figure~\ref{fig:lemma4-21-a1}, right side. 

Now assume that $a_1'|a_2'a_3'$; then $a_3'|a_4'a_5'$ and $a_1'|a_2'a_5'$ by canonicity. The latter implies $a_1'|a_3'a_4'$, and thus $a_2'|a_3'a_4'$ again by canonicity.  This leaves us with the following possibility:
\begin{itemize}
\item[(A2)] $a_1'|a_2'a_5'$,  $a_2'|a_3'a_5'$, $a_3'|a_4'a_5'$.
\end{itemize}

Next let $b_1,\ldots,b_5\in \mS_2$ be so that $b_1\prec\cdots\prec b_5$ and so that $b_1\perp b_4$, $b_2,b_3<b_4$, $b_2\perp b_3$, and  $b_1,\ldots,b_4<b_5$; see Figure~\ref{fig:lemma4-21-b1}, left side. 
\begin{figure}
\begin{center}
\includegraphics[scale=.5]{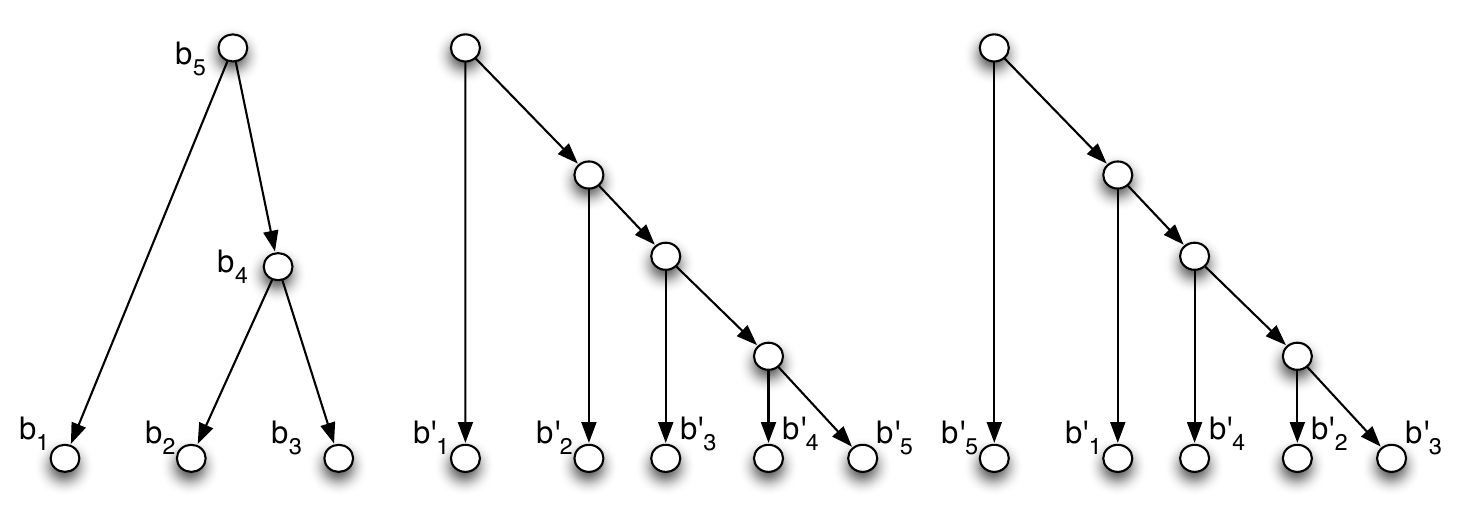}
\end{center}
\caption{Second illustration for the analysis of the behaviour of the function $f$ from Lemma~\ref{lem:flatemb}; the picture in the middle corresponds to case (B1), the picture on the right to case (B4).}
\label{fig:lemma4-21-b1}
\end{figure}

If $b_2'|b_3'b_4'$, then canonicity implies $b_1'|b_2'b_5'$ and $b_2'|b_3'b_5'$ leaving us with only two possibilities, namely $b_3'|b_4'b_5'$ and $b_3'b_4'|b_5'$. 
\begin{itemize}
\item[(B1)] $b_1'|b_2'b_5'$, $b_2'|b_3'b_5'$, $b_3'|b_4'b_5'$;
\item[(B2)] $b_1'|b_2'b_5'$, $b_2'|b_3'b_5'$, $b_3'b_4'|b_5'$.
\end{itemize}
See Figure~\ref{fig:lemma4-21-b1}, right side, for an illustration of case (B1). 
If on the other hand $b_2'b_3'|b_4'$, then canonicity tells us that $b_1'b_4'|b_5'$. One possibility here is that $b_1'b_2'|b_3'$, which together with $b_2'b_3'|b_4'$ implies $b_1'b_3'|b_4'$, and so we have:
\begin{itemize}
\item[(B3)] $b_1'b_4'|b_5'$, $b_1'b_3'|b_4'$, $b_1'b_2'|b_3'$.
\end{itemize}
Finally, suppose that $b_2'b_3'|b_4'$ and $b_1'|b_2'b_3'$. Pick $x>b_3$ such that $b_2\perp x$. Then $b_1'|b_2'x'$ by canonicity, and hence $b_2'\prec b_3'\prec x$ implies that we must have $b_1'|b_3'x'$. But then canonicity gives us $b_1'|b_2'b_4'$, and hence the following:
\begin{itemize}
\item[(B4)] $b_1'b_4'|b_5'$, $b_1'|b_2'b_4'$, $b_2'b_3'|b_4'$.
\end{itemize}

We now consider all possible combinations of these situations. Assume first that (A1) holds; then neither (B1) nor (B2) hold because otherwise $a_1'a_4'|a_5'$ and $b_1'|b_4'b_5'$ together would contradict canonicity. If we have (B3), then for all $a,b,c$ in the range of $f$ we have that $ab|c$ \iff $a,b\prec c$. 
We claim that then
$\Gamma$ is a reduct of $(\st;\prec)$. To see this, recall that every relation of $\Gamma$
has a quantifier-free definition over $(\st;\leq,C)$ and hence also a quantifier-free 
definition $\phi$ 
over $(\st;<,C)$. Note that when we evaluate 
$\phi$ on elements from the image of the flat function $f$, 
then atomic formulas of the form $x < y$ in $\phi$  
can be replaced by `false' without changing
the truth value of the formula. 
 Since $f$
is an embedding and preserves $\phi$, 
this shows that 
we can assume that $\phi$
does not make use of $<$. 
Moreover, since $f$ preserves $\prec$, 
we can replace occurrences
of $C(c,ab)$ in $\phi$ by $a \prec c \wedge b \prec c$, and obtain a formula which defines the same relation over $(\st;\prec)$, 
and this proves the claim. 
The structure $(\st;\prec)$ is isomorphic to $(\mathbb Q;<)$, and it follows
that $\Gamma$ is isomorphic to a reduct of $(\mathbb Q;<)$.
If we have (B4), then $f$ is a flat self-embedding of $\Gamma$, and $f$ preserves $R$:
to see this, let $(x,y,z) \in R$.
If $z|xy$ the $z'|x'y'$; if $x<z \wedge y < z$ 
then $z'|x'y'$; otherwise, $x \perp z \wedge y \perp z$
and $x$ and $y$ are comparable, 
and we again have $z'|x'y'$. 
In all cases, $(x',y',z') \in R$. 

Now assume that (A2) holds. Then $a_1'|a_4'a_5'$ and canonicity imply that (B1) or (B2) is the case. However, (B2) is in fact impossible by virtue of $a_1'a_3'|a_5'$ and $b_2'|b_4'b_5'$, leaving us with (B1). Here, we argue that $\Gamma$ is isomorphic to a reduct of $(\mathbb Q;<)$ precisely as in the case (A1)+(B3).
\end{proof}

\begin{lemma}\label{lem:funny}
Let $\Gamma$ be a reduct of $(\mS_2; \leq)$. Assume that there is a flat function in $\overline{\Aut(\Gamma)}$ that preserves $R$. Then $\Gamma$ is isomorphic to a reduct of $(\mathbb Q;<)$.
\end{lemma}
\begin{proof}
Let $f$ be that function. We use induction to show that the action 
of $\Aut(\Gamma)$ is 
\emph{highly set-transitive}, 
i.e., if two subsets of $\st$ have the same finite cardinality $n$, then there exists an automorphism of $\Gamma$ sending one set to the other; in other words, the setwise action of $\Aut(\Gamma)$ on $n$-element subsets of $\mS_2$ is transitive, for all $n \geq 1$. The statement is obvious for $n=1$ since
already $(\mS_2;\leq)$ is transitive. 
For $n=2$, the same argument does not work: 
the structure $(\mS_2;\leq)$
has two orbits of $n$-element subsets, and by the homogeneity of $(\mS_2;\leq,C)$,
this is the orbit of two comparable elements and the orbit of two incomparable elements in $(\mS_2;\leq)$.
However, the claim is also easy to see, because $f$ maps comparable elements $u,v$ in $(\mS_2;\leq)$ to incomparable elements, and 
$f \in \overline{\Aut(\Gamma)}$:
hence, $\{u,v\}$ lies in the same orbit as $\{f(u),f(v)\}$. 

Assume that the claim holds for some $n\in \mathbb{N}$, and let $A_1, A_2$ be $(n+1)$-element subsets and $a_i\in A_i$ for $i\in \{1,2\}$. By the induction hypothesis, for all $i\in \{1,2\}$ there exists an $\alpha_i\in \Aut(\Gamma)$ such that $\alpha_i[A_i\setminus \{a_i\}]$ is a chain. 
Hence we can extend the order $\leq$ on $A_i$ 
to a linear order $\sqsubseteq_i$ such
that $y \sqsubseteq a_i$ 
for all $y \in A_i$ with $y \perp a_i$. 
Then $$C(f \circ \alpha_i (u),f \circ \alpha_i(v) 
f \circ \alpha_i (w)) \Leftrightarrow \big (
f \circ \alpha_i (v)  \sqsubseteq_i f \circ 
\alpha_i (u) \text{ and } f \circ \alpha_i (w)  \sqsubseteq_i f \circ
\alpha_i (u) \big )$$ 
because $f$ preserves $R$. 
So $C$ on $f \alpha_i[A_i]$ is completely determined by $\sqsubseteq_i$, and in the same way
for $i=1$ and $2$, so $f \alpha_1[A_1]$ and $f \alpha_2[A_2]$ induce isomorphic structures in $(\mS_2;\leq,C)$.
The homogeneity of $(\mS_2;\leq,C)$ then implies
that there exists $\gamma \in 
\Aut(\mS_2; \leq, C) \subseteq \Aut(\Gamma)$ such that $\gamma[(f\circ \alpha_1)[A_1]]=(f\circ \alpha_2)[A_2]$. Hence, $\beta_2^{-1}\circ \gamma\circ \beta_1[A_1]= A_2$.
As $\Gamma$ is highly set-transitive, the assertion follows from Cameron's theorem~\cite{Cameron5} 
which states that every highly set-transitive structure is isomorphic to
a reduct of $({\mathbb Q};<)$. 
\end{proof}

\begin{proof}[Proof of Theorem~\ref{thm:groups}]
Let $\Gamma'$ be the structure that we obtain from $\Gamma$ by adding all first-order definable relations in $\Gamma$. Then $\Aut(\Gamma)= \Aut(\Gamma')$ and  $\overline{ \Aut(\Gamma')}= \Emb(\Gamma')=\End(\Gamma')$ since $\Gamma'$ is a model-complete core. 
Hence, if $\End(\Gamma')$ contains a flat 
function, then $\Gamma'$ is isomorphic to a reduct of $(\mQ;<)$ by Lemma~\ref{lem:flatemb} and Lemma~\ref{lem:funny}, and we are done. 
If $\End(\Gamma')$ contains a thin
function, then $\Gamma'$ is isomorphic to a reduct of $(\mQ;\leq)$ by Lemma~\ref{lem:thinemb}. 
Otherwise, by 
Theorem~\ref{thm:4case}
we have that $\End(\Gamma') \in \{ \overline{\Aut(\mS_2; \leq)}, \overline{\Aut(\mS_2; B)}\}$,
and thus $\Aut(\Gamma) = \Aut(\Gamma') \in \{ \Aut(\mS_2; \leq), \Aut(\mS_2; B)\}$
which concludes the proof of the statement. 
%
\end{proof}

\section{Applications in Constraint Satisfaction}\label{sect:csp}
Let $\Gamma$ be a structure with a finite relational signature $\tau$.
Then $\Csp(\Gamma)$, the \emph{constraint satisfaction problem} for $\Gamma$,
is the computational problem of deciding for a given finite $\tau$-structure
whether there exists a homomorphism to $\Gamma$.
There are several computational problems in the literature that
can be formulated as CSPs for reducts of $(\mS_2;\leq)$. 

When $\Gamma_b$ is the reduct of $(\mS_2;\leq)$ that contains precisely the \emph{binary} relations with a first-order definition in $(\mS_2;\leq)$,
then $\Csp(\Gamma_b)$ has been studied
under the name ``network consistency problem for the 
branching-time relation algebra'' by Hirsch~\cite{Hirsch}; it is shown
there that the problem can be solved in polynomial time.
For concreteness, we mention 
that in particular the problem $\Csp(\mS_2;<,\perp)$ 
can be solved in polynomial time, since it can be seen as a special case of $\Csp(\Gamma_b)$. 
Broxvall and Jonsson~\cite{BroxvallJonsson} found a better
algorithm for $\Csp(\Gamma_b)$ which improves the running
time from $O(n^5)$ to $O(n^{3.326})$, where $n$ is the number of elements in the input structure. 
Yet another algorithm
with a running time that is quadratic in the input size has
been described in~\cite{BodirskyKutz}. 
The complexity of the CSP of \emph{disjunctive reducts} of $(\mS_2;\leq,\prec)$ 
has been determined in~\cite{BroxvallJonsson}; a \emph{disjunctive
reduct} is a reduct each of whose relations can be defined by a disjunction
of the basic relations in such a way that the disjuncts do not share common variables. 

Independently from this line of research, motivated by research 
in computational linguistics, Cornell~\cite{Cornell}
studied the reduct $\Gamma_c$ of $(\mS_2;\leq,\prec)$ containing
all binary relations that are first-order definable over $(\mS_2;\leq,\prec)$.
Contrary to a conjecture of Cornell, it has been shown that 
$\Csp(\Gamma_c)$ (and in fact already $\Csp(\mS_2;<,\perp)$) cannot be solved by establishing path consistency~\cite{BodirskyKutzAI},
which is one of the most studied  algorithmic approaches to solve infinite-domain CSPs. However, $\Csp(\Gamma_c)$ can be solved in polynomial time~\cite{BodirskyKutzAI}. 

It is a natural but challenging research question 
to ask for a classification of the complexity
of $\Csp(\Gamma)$ for all reducts of $(\mS_2;\leq)$. In this context, we 
call the reducts of $(\mS_2;\leq)$ \emph{tree description constraint languages}. 
Such classifications have been obtained for the reducts of $({\mathbb Q};\leq)$
and the reducts of the random graph~\cite{tcsps-journal,BodPin-Schaefer-both}. 
In both these previous classifications, the classification of the model-complete cores
of the reducts played a central role. Our 
Theorem~\ref{thm:4case} shows that every tree description language 
belongs to at least one out of four cases; in cases one and two, the CSP
has already been classified. It is easy to show (and this will appear in forthcoming work) that the CSP is NP-hard when case three 
of Theorem~\ref{thm:4case} applies. It is also easy to see
(again we have to refer to forthcoming work) that in case four
of Theorem~\ref{thm:4case}, adding the relations $<$ and $\perp$
to $\Gamma$ does not change the computational complexity of the CSP. The corresponding fact for the reducts of $({\mathbb Q};\leq)$
and the reducts of the random graph has been extremely useful
in the subsequent classification. 
Therefore, the present paper and in particular 
Theorem~\ref{thm:4case} are highly relevant for the study of the CSP for tree description constraint languages.

\section*{Acknowledgements}
We would like to thank Dietrich Kuske for spotting a mistake in an earlier version of the article, Manfred Droste for valuable hints concerning the literature on semilinear orders, 
Miodrag Soki\'c and Jan Hubi\v{c}ka for helpful pointers to relevant literature on structural Ramsey properties, 
and the anonymous referee for comments that helped us to improve the presentation. 

\bibliographystyle{alpha}
\bibliography{local.bib}

\end{document}